\documentclass{icmart}

\contact[Anton.Zorich@univ-rennes1.fr]{IRMAR, Universit\'e
Rennes 1, Campus de Beaulieu, 35042 Rennes, FRANCE}

\newcommand{\Closure}{{\mathcal N}}
\newcommand{\noz}{m}
\newcommand{\noi}{n}
\newcommand{\Screen}{W}
\newcommand{\Subspace}{L}
\newcommand{\Cocycle}{A}
\newcommand{\SVc}{c(\cC)}
\renewcommand{\Re}{\operatorname{Re}}

\newcommand{\Vol}{\operatorname{Vol}}

\newcommand{\dist}{\operatorname{dist}}

\newcommand{\ind}{\operatorname{ind}}
\newcommand{\id}{\operatorname{id}}
\newcommand{\Id}{\operatorname{Id}}

\newcommand{\SO}{\operatorname{SO}(2,\mathbb R)}
\newcommand{\SLZ}{\operatorname{SL}(2,\mathbb Z)}
\newcommand{\SLR}{\operatorname{SL}(2,\mathbb R)}
\newcommand{\GL}{\operatorname{GL}}
\newcommand{\SL}{\operatorname{SL}}

\newcommand\Z[1]{{\mathbb Z}^{#1}}

\newcommand\R[1]{{\mathbb R}^{#1}}
\newcommand\C[1]{{\mathbb C}^{#1}}
\newcommand\torus{{\mathbb T}^2}

\newcommand\Hyp{{\mathbb H}^2}

\newcommand{\cC}{{\mathcal C}}
\newcommand{\cG}{{\mathcal G}}
\newcommand{\cH}{{\mathcal H}}
\newcommand{\cM}{{\mathcal M}}
\newcommand{\cQ}{{\mathcal Q}}
\newcommand{\cS}{{\mathcal S}}
\newcommand{\cT}{{\mathcal T}}

\newlength{\halfbls}

\newtheorem{Theorem}{Theorem}
\newtheorem*{NNTheorem}{Theorem}
\newtheorem*{KeyTheorem}{Key~Theorem}
\newtheorem{TheoremPrime}{Theorem}

\newtheorem*{NNLemma}{Lemma}

\theoremstyle{definition}
\newtheorem*{NNRemark}{Remark}

\title[Geodesics on Flat Surfaces]{Geodesics on Flat Surfaces}

\author[Anton Zorich]{Anton Zorich\thanks{The author is grateful to
the MPI (Bonn) and to IHES for hospitality during the preparation of
this paper.}}

\begin{document}

\begin{picture}(0,0)(0,0)
\put(170,-490){\small Proceedings of the International Congress}
\put(170,-500){\small of Mathematics, Madrid, Spain, 2006}
\put(170,-510){\small 2006, EMS Publishing House}
\end{picture}

\begin{abstract}
Various problems of  geometry,  topology and dynamical systems on
surfaces as  well  as  some  questions concerning one-dimensional
dynamical systems lead to the  study  of  closed surfaces endowed
with a flat metric with  several  cone-type  singularities. In an
important  particular  case,  when  the flat metric  has  trivial
holonomy, the corresponding flat surfaces are naturally organized
into families  which appear to  be isomorphic to moduli spaces of
holomorphic one-forms.

One can obtain much  information  about the geometry and dynamics
of an  individual flat surface by  studying both its  orbit under
the Teichm\"uller geodesic flow and under the linear group action
on  the  corresponding  moduli  space.  We   apply  this  general
principle to  the study of generic  geodesics and to  counting of
closed geodesics on a flat surface.
\end{abstract}



\maketitle

\section*{Introduction: Families of Flat Surfaces as
Moduli Spaces of Abelian Differentials}

\begin{figure}[htb]
%
%
\includegraphics{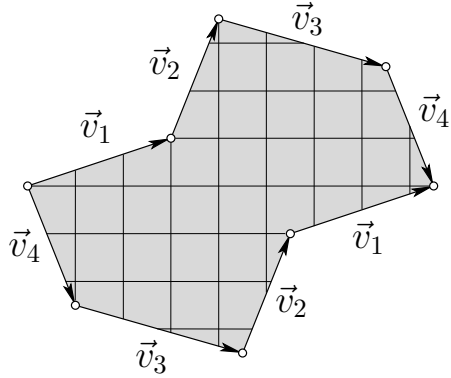}
\begin{picture}(20,2)(20,2)
\put(142,-48){\Large $\vec v_1$}
\put(168,-24){\Large $\vec v_2$}
\put(222,-7){\Large $\vec v_3$}
\put(270,-42){\Large $\vec v_4$}
\put(115,-93){\Large $\vec v_4$}
\put(163,-135){\Large $\vec v_3$}
\put(216,-114){\Large $\vec v_2$}
\put(245,-90){\Large $\vec v_1$}
\end{picture}
\vspace{140bp}
\caption{
\label{zorich:fig:suspension}
Identifying  corresponding  pairs  of  sides of this  polygon  by
parallel translations  we obtain a surface  of genus two.  It has
single  conical  singularity  with  cone angle $6\pi$;  the  flat
metric has trivial holonomy.
}
\end{figure}

Consider a collection of vectors $ \vec v_1,  \dots, \vec v_\noi$
in $\R{2}$ and construct from  these  vectors a broken line in  a
natural way:  a $j$-th edge of the broken  line is represented by
the vector $\vec{v}_j$. Construct another broken line starting at
the same point as  the initial one by taking the same  vectors in
the order  $\vec{v}_{\pi(1)}, \dots, \vec{v}_{\pi(\noi)}$,  where
$\pi$ is some permutation of $\noi$ elements. By construction the
two  broken  lines share the same endpoints;  suppose  that  they
bound   a   polygon   like in   Fig.~\ref{zorich:fig:suspension}.
Identifying the pairs  of sides corresponding to the same vectors
$\vec{v}_j$, $j=1,  \dots,  \noi$,  by  parallel  translations we
obtain a surface endowed with  a  flat  metric. (This construction
follows  the  one   in~\cite{zorich:Masur:Annals:82}.)  The  flat
metric is nonsingular outside  of  a finite number of cone-type
singularities  corresponding  to the vertices of the polygon.  By
construction the flat metric has  trivial  holonomy:  a  parallel
transport of a  vector  along a closed path  does  not change the
direction   (and   length)  of  the  vector.  This  implies,   in
particular, that all cone angles are    integer  multiples  of
$2\pi$.

The  polygon  in  our  construction depends continuously  on  the
vectors $\vec{v}_j$.  This  means that the combinatorial geometry
of the  resulting flat surface  (its genus $g$, the number $\noz$
and types of the resulting conical singularities) does not change
under small deformations of the vectors  $\vec{v}_j$. This allows
to consider  a flat surface  as an  element of a  family of  flat
surfaces  sharing  common combinatorial geometry; here we do  not
distinguish isometric flat surfaces. As an example of such family
one can consider a family of  flat tori of area one, which can be
identified with the space of lattices of area one:
$$
\begin{array}{c}
\backslash\, \SLR\, /\\
[-\halfbls]
\SO\ \phantom{\backslash\, \SLR\! /} \SLZ
\end{array}
\quad
=
\quad
\begin{array}{l}
\Hyp /\\
[-\halfbls]\phantom{\Hyp /}   \SLZ
\end{array}
$$

\begin{figure}[htb]
%
\centering
\includegraphics{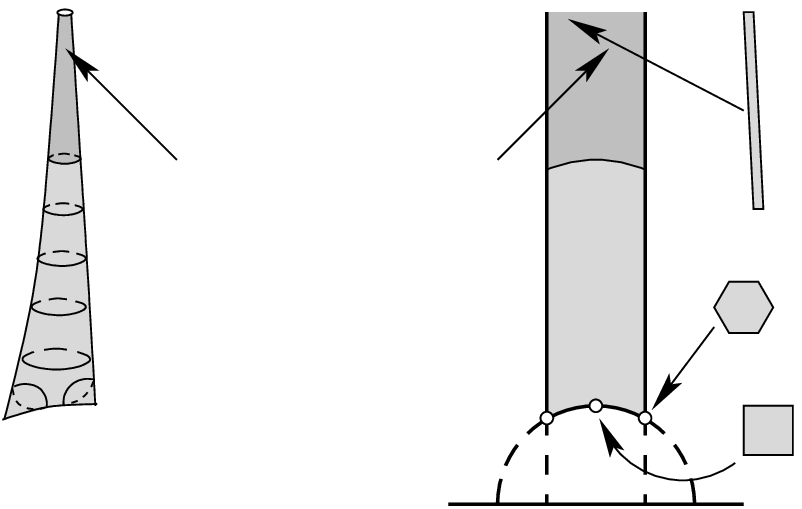}
\begin{picture}(-11,0)(-11,0)
\put(-70,-40){neighborhood of a}
\put(-70,-50){cusp = subset of}
\put(-70,-60){tori having short}
\put(-70,-70){closed geodesic}
\end{picture}
\vspace{140bp} 
\caption{
\label{zorich:fig:space:of:flat:tori}
``Modular surface''  $\Hyp/\SLZ$  representing the space of flat
tori is a noncompact orbifold of finite volume.
}
\end{figure}

The   corresponding   ``modular   surface''   is   not   compact,
see~Fig.~\ref{zorich:fig:space:of:flat:tori}.      Flat      tori
representing  points,  which are close to the  cusp,  are  almost
degenerate: they have a very short closed  geodesic. Similarly,
families of flat surfaces  of  higher genera also form noncompact
finite-dimensional orbifolds. The origin of their  noncompactness
is the same as for  the  tori: flat surfaces having short  closed
geodesics   represent    points    which   are   close   to   the
multidimensional ``cusps''.

We shall consider  only those flat  surfaces,  which
have trivial holonomy. Choosing a direction at some point of such
flat surface we can transport it  to any other point. It would be
convenient to include the  choice  of direction in the definition
of a  flat structure. In  particular, we want to distinguish the flat
structure            represented by the            polygon
in~Fig.~\ref{zorich:fig:suspension} and the one represented by the
same polygon rotated by  some  angle different from $2\pi$.

Consider the natural coordinate $z$ in the complex plane. In this
coordinate the parallel translations which we use to identify the
sides of the polygon in~Fig.~\ref{zorich:fig:suspension} are represented as $z'=z+const$.
Since this correspondence  is holomorphic, it means that our flat
surface  $S$  with punctured conical points inherits the  complex
structure. It is easy to check that the complex structure extends
to the punctured points.  Consider  now a holomorphic 1-form $dz$
in  the complex plane.  When  we  pass  to the  surface  $S$  the
coordinate  $z$  is not globally defined anymore. However,  since
the changes of  local coordinates are defined as $z'=z+const$, we
see that $dz=dz'$. Thus, the  holomorphic  1-form  $dz$ on $\C{}$
defines  a  holomorphic  1-form  $\omega$ on $S$ which  in  local
coordinates has  the form $\omega=dz$. It  is easy to  check that
the form $\omega$ has zeroes exactly at those points of $S$ where
the flat structure has conical singularities.

Reciprocally, one can
show  that  a pair (Riemann surface, holomorphic  1-form)
uniquely defines a flat structure of the type described above.

In an appropriate local  coordinate  $w$ a holomorphic 1-form can
be represented in a neighborhood of zero as  $w^d\,dw$, where $d$
is called  the degree of zero. The  form $\omega$  has a zero  of
degree $d$  at a conical point  with cone angle  $2\pi(d+1)$. The
sum of  degrees  $d_1+\dots+d_\noz$  of  zeroes  of a holomorphic
1-form  on a  Riemann  surface of genus  $g$  equals $2g-2$.  The
moduli space  $\cH_g$  of  pairs  (complex structure, holomorphic
1-form) is a  $\C{g}$-vector bundle over the moduli space $\cM_g$
of complex structures. The space $\cH_g$  is naturally stratified
by  the  strata  $\cH(d_1,\dots,d_\noz)$ enumerated by  unordered
partitions  of  the  number  $2g-2$ in a collection  of  positive
integers   $2g-2=d_1+\dots+d_\noz$.   Any   holomorphic   1-forms
corresponding  to  a fixed  stratum  $\cH(d_1,\dots,d_\noz)$  has
exactly $\noz$ zeroes, and $d_1,\dots,d_\noz$ are  the degrees of
zeroes. Note, that an individual stratum  $\cH(d_1,\dots,d_\noz)$
in general does not form a fiber bundle over $\cM_g$.

It is possible to show  that  if the permutation $\pi$ which  was
used to  construct a polygon  in~Fig.~\ref{zorich:fig:suspension}
satisfy  some  explicit conditions,  vectors  $\vec{v}_1,  \dots,
\vec{v}_\noi$ representing  the  sides  of  the  polygon serve as
coordinates in the  corresponding family $\cH(d_1,\dots,d_\noz)$.
Consider vectors $\vec v_j$ as  complex  numbers.  Let $\vec v_j$
join  vertices  $P_j$  and  $P_{j+1}$ of the polygon.  Denote  by
$\rho_j$  the  resulting  path  on   $S$   joining   the   points
$P_j,P_{j+1}\in  S$.  Our interpretation of $\vec v_j$  as  of  a
complex number implies that
$$
\int_{\rho_j}\omega
=\int_{P_j}^{P_{j+1}}dz=
v_j \in \C{}
$$
The path $\rho_j$ represents a relative cycle: an  element of the
relative homology group $H_1(S,\{P_1,\dots,P_\noz\}\,;\,\Z{})$ of
the  surface  $S$ relative to the finite  collection  of  conical
points $\{P_1,\dots,P_\noz\}$.  Relation  above  means that $\vec
v_j$ represents a period of $\omega$:  an  integral  of  $\omega$
over the relative cycle $\rho_j$. In other words, a small domain in
$H^1(S,\{P_1,\dots,P_\noz\};\C{})$ containing $[\omega]$  can  be
considered  as   a   local   coordinate   chart   in  our  family
$\cH(d_1,\dots,d_\noz)$ of flat surfaces.

We summarize  the  correspondence  between  geometric language of
flat surfaces  and  the  complex-analytic language of holomorphic
1-forms on a Riemann surface in the dictionary below.

\bigskip

\begin{tabular}{|c|c|}
\hline&\\[-\halfbls]
       Geometric language              &      Complex-analytic language\\
[-\halfbls] &\\ \hline & \\ [-\halfbls]
flat structure (including a choice     &      complex structure and a choice\\
of the vertical direction)             &   of a holomorphic 1-form $\omega$\\
[-\halfbls] &\\ \hline & \\ [-\halfbls]
conical point                          & zero of degree $d$\\
with a cone angle $2\pi(d+1)$          & of the holomorphic 1-form $\omega$\\
                                       &(in local coordinates $\omega=w^d\,dw$)\\
[-\halfbls] &\\ \hline & \\ [-\halfbls]
side $\vec v_j$ of a polygon           & relative period $\int_{P_j}^{P_{j+1}}\omega=\int_{\vec v_j}\omega$ \\
                                       & of the 1-form $\omega$ \\
[-\halfbls] &\\ \hline & \\ [-\halfbls]
family of flat surfaces sharing        & stratum $\cH(d_1,\dots,d_\noz)$ in the \\
the same cone angles                   &  moduli space of Abelian differentials \\
$2\pi(d_1+1),\dots,2\pi(d_\noz+1)$     & \\
[-\halfbls] &\\ \hline & \\ [-\halfbls]
coordinates in the family:             & coordinates in $\cH(d_1,\dots,d_\noz)$ : \\
vectors $\vec v_i$                     & relative periods of $\omega$ in  \\
defining the polygon                   & $H^1(S, \{P_1,\dots,P_\noz\};\C{})$ \\
[-\halfbls] &\\ \hline
\end{tabular}

\bigskip

Note    that    the    vector    space    $H^1(S,\{P_1,    \dots,
P_\noz\}\,;\,\C{})$   contains   a    natural   integer   lattice
$H^1(S,\{P_1,  \dots,  P_\noz\}\,;\,\Z{}\oplus \sqrt{-1}\,\Z{})$.
Consider  a  linear  volume  element $d\nu$ in the  vector  space
$H^1(S,\{P_1, \dots, P_\noz\}\,;\,\C{})$ normalized in such a way
that  the  volume  of  the  fundamental  domain in the  ``cubic''
lattice
$$
H^1(S,\{P_1, \dots, P_\noz\}\,;\,\Z{}\oplus \sqrt{-1}\,\Z{})\ \subset\
H^1(S,\{P_1, \dots, P_\noz\}\,;\,\C{})
$$
equals    one.    Consider    now    the    real     hypersurface
$\cH_1(d_1,\dots,d_\noz)\subset\cH(d_1,\dots,d_\noz)$ defined  by
the  equation  $area(S)=1$. The  volume  element  $d\nu$  can  be
naturally restricted  to  the  hypersurface  defining  the volume
element $d\nu_1$ on $\cH_1(d_1,\dots,d_\noz)$.

\begin{NNTheorem}[H.~Masur. W.~A.~Veech]
The total volume $\Vol(\cH_1(d_1,\dots,d_\noz))$ of every stratum
is finite.
\end{NNTheorem}
The values of these volumes  were  computed by A.~Eskin
and A.~Okounkov~\cite{zorich:Eskin:Okounkov}.

Consider  a  flat surface $S$ and consider  a  polygonal  pattern
obtained by unwrapping $S$  along some geodesic cuts. For example,
one can assume that our flat  surface $S$ is glued from a polygon
$\Pi\subset   \R{2}$   as  on   Fig.~\ref{zorich:fig:suspension}.
Consider a linear transformation $g\in \GL^+(2,\R{})$ of the plane
$\R{2}$. The sides of the new polygon $g \Pi$ are  again arranged
into pairs, where the sides  in  each pair are parallel and  have
equal length.  Identifying the sides in  each pair by  a parallel
translation we  obtain a new  flat surface $g S$ which, actually,
does  not depend  on the  way  in which  $S$ was  unwrapped to  a
polygonal pattern $\Pi$. Thus, we get a continuous  action of the
group $\GL^+(2,\R{})$ on each stratum $\cH(d_1,\dots,d_\noz)$.

Considering the  subgroup  $\SL(2,\R{})$ of area preserving linear
transformations we get the  action  of $\SL(2,\R{})$ on the ``unit
hyperboloid'' $\cH_1(d_1,\dots,d_\noz)$. Considering the diagonal
subgroup $\begin{pmatrix}  e^t  &  0\\  0  & e^{-t} \end{pmatrix}
\subset   \SL(2,\R{})$   we  get  a  continuous  action  of   this
one-parameter subgroup  on each stratum  $\cH(d_1,\dots,d_\noz)$.
This  action induces  a  natural flow on  the  stratum which  is
called the \emph{Teichm\"uller geodesic flow}.

\begin{KeyTheorem}[H.~Masur. W.~A.~Veech]
The actions of  the groups $\SL(2,\R{})$ and $\begin{pmatrix} e^t &
0\\ 0 & e^{-t} \end{pmatrix}$  preserve  the  measure  $d\nu_1$.
Both actions are  ergodic with respect  to this measure  on  each
connected component of every stratum $\cH_1(d_1,\dots,d_\noz)$.
\end{KeyTheorem}

The following  basic principle (which was  was first used  in the
pioneering works of H.~Masur~\cite{zorich:Masur:Annals:82} and of
W.~Veech~\cite{zorich:Veech:Annals:82} to prove unique ergodicity
of almost all  interval  exchange transformations) appeared to be
surprisingly powerful in the study of flat surfaces. Suppose that
we  need  some  information  about  geometry  or  dynamics  of an
individual  flat  surface $S$. Consider the ``point''  $S$  in  the
corresponding family  of flat surfaces $\cH(d_1, \dots, d_\noz)$.
Denote  by  $\Closure(S)=\overline{\GL^+(2,\R{})\,  S}\subset  \cH(d_1,
\dots, d_\noz)$ the closure of the $\GL^+(2,\R{})$-orbit of $S$ in
$\cH(d_1,  \dots,  d_\noz)$.

In  numerous  cases  knowledge
about the structure of $\Closure(S)$ gives a comprehensive information
about  geometry  and dynamics of the initial  flat  surface  $S$.
Moreover, some delicate  numerical  characteristics of $S$ can be
expressed as  averages of simpler characteristics over $\Closure(S)$.
We apply this  general philosophy to  the study of  geodesics  on
flat surfaces.

Actually,  there  is  a  hope that this philosophy  extends  much
further. A  closure of an orbit  of an abstract  dynamical system
might have  extremely  complicated  structure.  According  to the
optimistic hopes,  the closure $\Closure(S)$ of a $\GL^+(2,\R{})$-orbit
of any  flat surface $S$  is a nice complex-analytic variety, and
all such varieties might be classified. For genus  two the latter
statements     were    recently     proved     by     C.~McMullen
(see~\cite{zorich:McMullen:Hilbert}
and~\cite{zorich:McMullen:genus:2})      and      partly       by
K.~Calta~\cite{zorich:Calta}.

The following theorem supports the hope for some  nice and simple
description of orbit closures.

\begin{NNTheorem}[M.~Kontsevich]
Suppose that the closure in the stratum $\cH(d_1, \dots, d_\noz)$ of a
$\GL^+(2,\R{})$-orbit   of    some   flat   surface   $S$   is   a
complex-analytic  subvariety.  Then in  cohomological coordinates
$H^1(S,\{P_1,\dots,P_\noz\};\C{})$ this subvariety is represented by an affine
subspace.
\end{NNTheorem}

\section{Geodesics Winding up Flat Surfaces}
\label{zorich:s:Generic:Geodesics}

In this section we study geodesics on a flat surface $S$ going in
generic directions.  According  to  the  theorem of S.~Kerckhoff,
H.~Masur  and   J.~Smillie~\cite{zorich:Kerckhoff:Masur:Smillie},
for  any flat  surface  $S$ the directional  flow  in almost  any
direction is uniquely ergodic. This implies,  in particular, that
for such directions the geodesics wind around $S$ in a relatively
regular manner. Namely,  it  is possible  to  find a cycle  $c\in
H_1(S;\R{})$ such that a long piece of geodesic  pretends to wind
around  $S$  repeatedly  following  this  asymptotic  cycle  $c$.
Rigorously  it  can be described as follows.  Having  a  geodesic
segment $X\subset S$ and  some point $x\in X$ we emit from  $x$ a
geodesic transversal to $X$. From time to time the geodesic would
intersect $X$.  Denote  the  corresponding  points  as $x_1, x_2,
\dots $.  Closing up the  corresponding pieces of the geodesic by
joining the  starting point $x_0$ and  the point $x_j$  of $j$-th
return  to  $X$  with a path going along $X$ we get a sequence of
closed paths defining  the cycles $c_1, c_2, \dots$. These cycles
represent  longer  and longer pieces of the  geodesic.  When  the
direction of the geodesic is uniquely ergodic, the limit
$$
\lim_{N\to\infty} \frac{1}{N}\ c_N =  c
$$
exists and  the corresponding asymptotic cycle $c\in H_1(S;\R{})$
does not  depend on the starting  point $x_0\in X$.  Changing the
transverse interval $X$ we get a collinear asymptotic cycle.

When $S$ is a flat torus glued from a unit square, the asymptotic
cycle $c$ is a  vector  in $H_1(\torus;\R{})=\R{2}$ and its
slope  is exactly  the  slope of our  flat  geodesic in  standard
coordinates. When $S$ is a surface of higher genus the asymptotic
cycle belongs  to a $2g$-dimensional space $H_1(S;\R{})= \R{2g}$.
Let us  study how the cycles $c_j$ deviate  from the direction of
the  asymptotic  cycle $c$.  Choose  a  hyperplane  $\Screen$  in
$H_1(S,\R{})$ orthogonal  (transversal)  to  the asymptotic cycle
$c$ and consider a parallel projection to  this screen along $c$.
Projections of the cycles $c_N$ would not be necessarily bounded:
directions  of  the  cycles
$c_N$ tend to direction of the asymptotic cycle  $c$ provided the
norms of the projections grow sublinearly with respect to $N$.

\begin{figure}[hbt]
\special
{
psfile=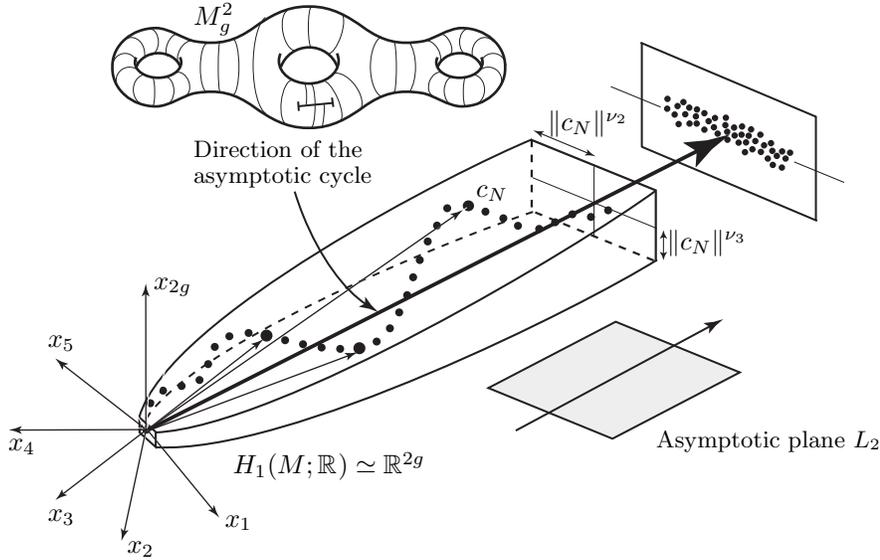
hscale= 100
vscale= 100
hoffset=20
voffset=-220
}
\begin{picture}(0,0)(0,0)
\put(0,0) 
 {\begin{picture}(0,0)(0,0)
 \put(197,-85){$c_N$}
 \put(224,-60){$\|c_N\|^{\nu_2}$}
 \put(269,-105){$\|c_N\|^{\nu_3}$}
\end{picture}}
\put(0,0)
 {\begin{picture}(0,0)(0,0)
   \put(105,-190){$H_1(M;\R{})\simeq\R{2g}$}
   \put(102,-210){$x_1$}
   \put(65,-220){$x_2$}
   \put(35,-208){$x_3$}
   \put(20,-181){$x_4$}
   \put(35,-142){$x_5$}
   \put(75,-120){$x_{2g}$}
 \end{picture}}
\put(0,0)
 {\begin{picture}(0,0)(0,0)
   \put(265,-180){\small Asymptotic plane $L_2$}
 \end{picture}}
\put(0,0)
 {\begin{picture}(0,0)(0,0)
   \put(90,-70){\small Direction of the}
   \put(90,-80){\small asymptotic cycle}
 \end{picture}}
\put(0,0)
 {\begin{picture}(0,0)(0,0)
   \put(90,-20){$M^2_g$}
 \end{picture}}
\end{picture}
\vspace{220bp}
\caption{
\label{zorich:fig:asymptotic:cycle:deviation}
Deviation  from  the asymptotic  direction  exhibits  anisotropic
behavior: vectors deviate mainly along two-dimensional  subspace,
a bit more along three-dimensional subspace, etc. Their deviation
from a  Lagrangian  $g$-dimensional subspace is already uniformly
bounded.
}
\end{figure}

Let us observe how the projections are distributed  in the screen
$\Screen$.      A       heuristic      answer      is       given
by~Fig.~\ref{zorich:fig:asymptotic:cycle:deviation}. We see  that
the distribution of projections of the cycles $c_N$ in the screen
$\Screen$ is  anisotropic:  the projections accumulate along some
line. This means that in the original space  $\R{2g}$ the vectors
$c_N$ deviate from the asymptotic direction $\Subspace_1$ spanned
by  $c$  not arbitrarily but along some two-dimensional  subspace
$\Subspace_2$ containing $\Subspace_1$,               see
Fig.~\ref{zorich:fig:asymptotic:cycle:deviation}.
Moreover, measuring  the norms $\|proj(c_N)\|$ of the projections
we get
$$
\limsup_{N\to\infty} \frac{\log\|proj(c_N)\|}{\log N} = \nu_2<1
$$
Thus,   the   vector  $c_N$  is  located  approximately  in   the
two-dimensional plane  $\Subspace_2$,  and  the distance from its
endpoint to  the line $\Subspace_1$ in $\Subspace_2$ is
at most of the order $\|c_N\|^{\nu_2}$,           see
Fig.~\ref{zorich:fig:asymptotic:cycle:deviation}.

Consider now a new screen $\Screen_2\perp\Subspace_2$  orthogonal
to  the  plane $\Subspace_2$.  Now  the  screen  $\Screen_2$  has
codimension   two  in   $H_1(S,\R{})\simeq\R{2g}$.   Taking   the
projections  of  $c_N$  to  $\Screen_2$  along  $\Subspace_2$  we
eliminate   the    asymptotic   directions   $\Subspace_1$    and
$\Subspace_2$ and  we see how     the vectors $c_N$  deviate from
$\Subspace_2$. On  the  screen  $\Screen_2$  we  observe the same
picture  as  in Fig.~\ref{zorich:fig:asymptotic:cycle:deviation}:
the  projections   are  again  located  along  a  one-dimensional
subspace.

Coming back  to the ambient space $H_1(S,\R{})\simeq\R{2g}$, this
means that  in the first  term of approximation all vectors $c_N$
are  aligned  along the  one-dimensional  subspace  $\Subspace_1$
spanned  by  the  asymptotic   cycle.   In  the  second  term  of
approximation, they  can  deviate  from  $\Subspace_1$,  but  the
deviation   occurs   mostly   in  the  two-dimensional   subspace
$\Subspace_2$, and  has order $\|c_N\|^{\nu_2}$ where $\nu_2<1$. In
the third term of approximation we see that the vectors $c_N$ may
deviate  from  the plane $\Subspace_2$, but the deviation  occurs
mostly in a  three-dimensional  space $\Subspace_3$ and has order
$\|c_N\|^{\nu_3}$ where $\nu_3<\nu_2$.

Going on we  get further terms of approximation. However, getting
to a  subspace $\Subspace_g$ which  has half of the dimension of the
ambient space  we  see that, in   there is no  more
deviation from  $\Subspace_g$:  the  distance  from  any $c_N$ to
$\Subspace_g$ is uniformly bounded.
Note   that    the    intersection    form   endows   the   space
$H_1(S,\R{})\simeq\R{2g}$ with a natural symplectic structure. It
can  be  checked  that  the  resulting  $g$-dimensional  subspace
$\Subspace_g$ is a Lagrangian subspace for this symplectic form.

A rigorous  formulation  of  phenomena described heuristically in
Fig.~\ref{zorich:fig:asymptotic:cycle:deviation} is  given by the
theorem below.

By convention  we always consider  a flat surface together with a
choice of direction which  is  called the vertical direction, or,
sometimes, ``direction  to  the  North''.  Using  an  appropriate
homotethy  we  normalize the area of  $S$  to one, so that  $S\in
\cH_1(d_1,\dots,d_\noz)$.

We chose a point $x_0\in S$ and a horizontal segment  $X$ passing
through  $x_0$;  by  $|X|$  we  denote  the  length  of  $X$.
The interval  $X$  is chosen in  such  way, that  the
interval exchange transformation induced by the vertical flow has
the minimal  possible  number $\noi=2g+\noz-1$ of subintervals
under exchange. (Actually, almost any other choice of $X$ would also work.)
We consider a  geodesic ray $\gamma$  emitted from $x_0$  in  the
vertical  direction.  (If  $x_0$  is  a  saddle point, there  are
several outgoing vertical geodesic rays;  choose  any  of  them.)
Each time when $\gamma$ intersects $X$ we join the point $x_N$ of
intersection and the starting  point  $x_0$ along $X$ producing a
closed path. We denote the  homology  class  of the corresponding
loop by $c_N$.

Let $\omega$ be  the holomorphic 1-form representing $S$; let $g$
be   genus   of   $S$.   Choose   some    Euclidean   metric   in
$H_1(S;\R{})\simeq\R{2g}$ which would allow to measure a distance
from  a  vector   to  a  subspace.  Let  by  convention  $\log(0)
=-\infty$.

\begin{Theorem}
\label{zorich:th:mf}
For    almost    any   flat   surface   $S$   in   any    stratum
$\cH_1(d_1,\dots,d_\noz)$ there exists a flag of subspaces
$$
\Subspace_1\subset \Subspace_2\subset\dots\subset \Subspace_g\subset H_1(S;\R{})
$$
in the first homology group of the surface with the following
properties.

Choose any starting point $x_0\in X$ in the horizontal segment
$X$. Consider the corresponding sequence $c_1, c_2, \dots $ of
cycles.

--- The following limit exists
$$
|X| \lim_{N\to\infty} \frac{1}{N}\, c_N = c,
$$
where the  nonzero  asymptotic  cycle  $c\in  H_1(M^2_g;\R{})$ is
Poincar\'e     dual     to    the     cohomology     class     of
$\omega_0=\Re[\omega]$,   and    the   one-dimensional   subspace
$\Subspace_1=\langle c\rangle_{\R{}}$ is spanned by $c$.

--- For any $j=1,\dots, g-1$ one has
$$
\limsup_{N\to\infty} \frac{\log\dist(c_N,\Subspace_j) }{\log N}  =
\nu_{j+1}
$$
and
$$
\dist(c_N,\Subspace_g)  \le const,
$$
where the constant depends  only on $S$ and on the choice  of the
Euclidean structure in the homology space.

The numbers  $2,1+\nu_2,\dots,1+\nu_g$  are  the top $g$ Lyapunov
exponents of the Teichm\"uller geodesic flow on the corresponding
connected         component         of        the         stratum
$\mathcal{H}(d_1,\dots,d_\noz)$;  in  particular,   they  do  not
depend  on  the  individual  generic  flat  surface  $S$  in  the
connected component.
\end{Theorem}

It  should  be stressed, that the theorem  above  was  formulated
in~\cite{zorich:Zorich:How:do} as a  conditional statement: under
the conjecture that  $\nu_g>0$  there exist a Lagrangian subspace
$\Subspace_g$ such that the cycles are in a bounded distance from
$\Subspace_g$;  under   the   further  conjecture  that  all  the
exponents $\nu_j$, for $j=2,\dots,g$, are distinct,  there is a
complete  Lagrangian  flag (i.e. the dimensions of the  subspaces
$\Subspace_j$, where  $j=1,2,\dots,g$,  rise  each  time by one).
These     two     conjectures    were     later     proved     by
G.~Forni~\cite{zorich:Forni:02}    and    by     A.~Avila     and
M.~Viana~\cite{zorich:Avila:Viana} correspondingly.

Currently  there  are  no  methods of calculation  of  individual
Lyapunov exponents  $\nu_j$  (though  there  is some experimental
knowledge of their approximate  values). Nevertheless,
for any connected  component of any stratum (and, more generally,
for any  $\GL^+(2;\R{})$-invariant  suborbifold) it is possible to
evaluate   the    \textit{sum}    of   the   Lyapunov   exponents
$\nu_1+\dots+\nu_g$, where $g$  is the genus. The formula for this
sum was discovered  by M.~Kontsevich;
morally, it  is given in terms  of characteristic  numbers of
some     natural     vector     bundles    over    the     strata
$\cH(d_1,\dots,d_\noz)$, see~\cite{zorich:Kontsevich}.
Another  interpretation of this  formula
was found  by  G.~Forni~\cite{zorich:Forni:02};  see  also a very
nice  formalization of   these   results   in   the   survey   of
R.~Krikorian~\cite{zorich:Krikorian}.     For    some     special
$\GL^+(2;\R{})$-invariant  suborbifolds  the corresponding  vector
bundles might have equivariant subbundles, which provides additional
information on  corresponding subcollections of the Lyapunov
exponents,  or even gives their explicit values in some  cases,
like in the case of Teichm\"uller curves considered in the paper
of I.~Bouw  and M.~M\"oller~\cite{zorich:Bouw:Moeller}.

Theorem~\ref{zorich:th:mf}    illustrates    a   phenomenon    of
\textit{deviation   spectrum}.   It  was   proved   by   G.~Forni
in~\cite{zorich:Forni:02} that ergodic sums  of  smooth functions
on  an   interval   along   trajectories   of  interval  exchange
transformations, and  ergodic  integrals  of  smooth functions on
flat  surfaces  along  trajectories  of  directional  flows  have
deviation  spectrum   analogous   to   the   one   described   in
theorem~\ref{zorich:th:mf}. L.~Flaminio and  G.~Forni showed that
the same phenomenon can be observed for other parabolic dynamical
systems, for example,  for the horocycle flow on compact surfaces
of constant negative curvature~\cite{zorich:Flaminio:Forni:03}.

\paragraph{Idea of the proof: renormalization.}

\begin{figure}
%
\includegraphics{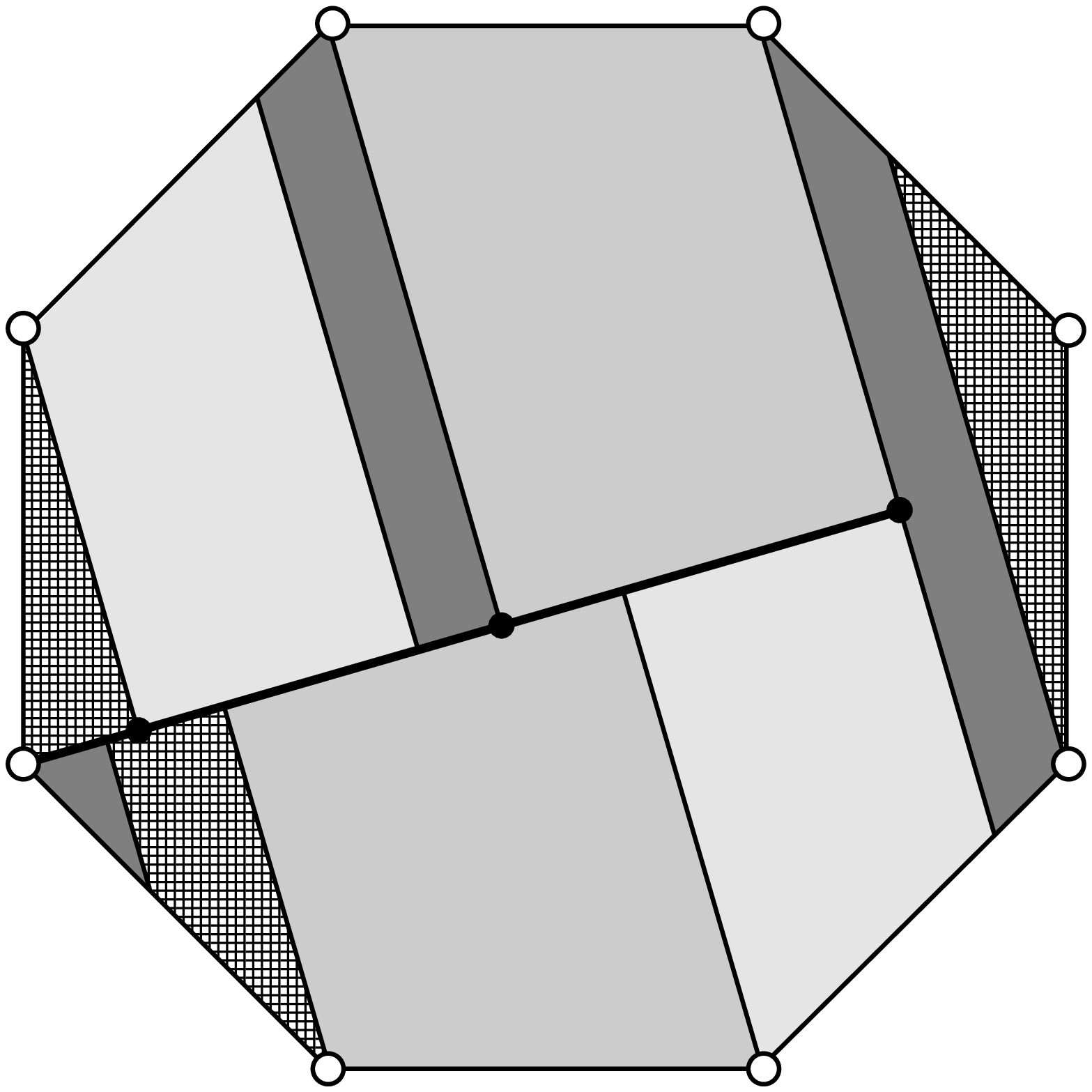}
\includegraphics{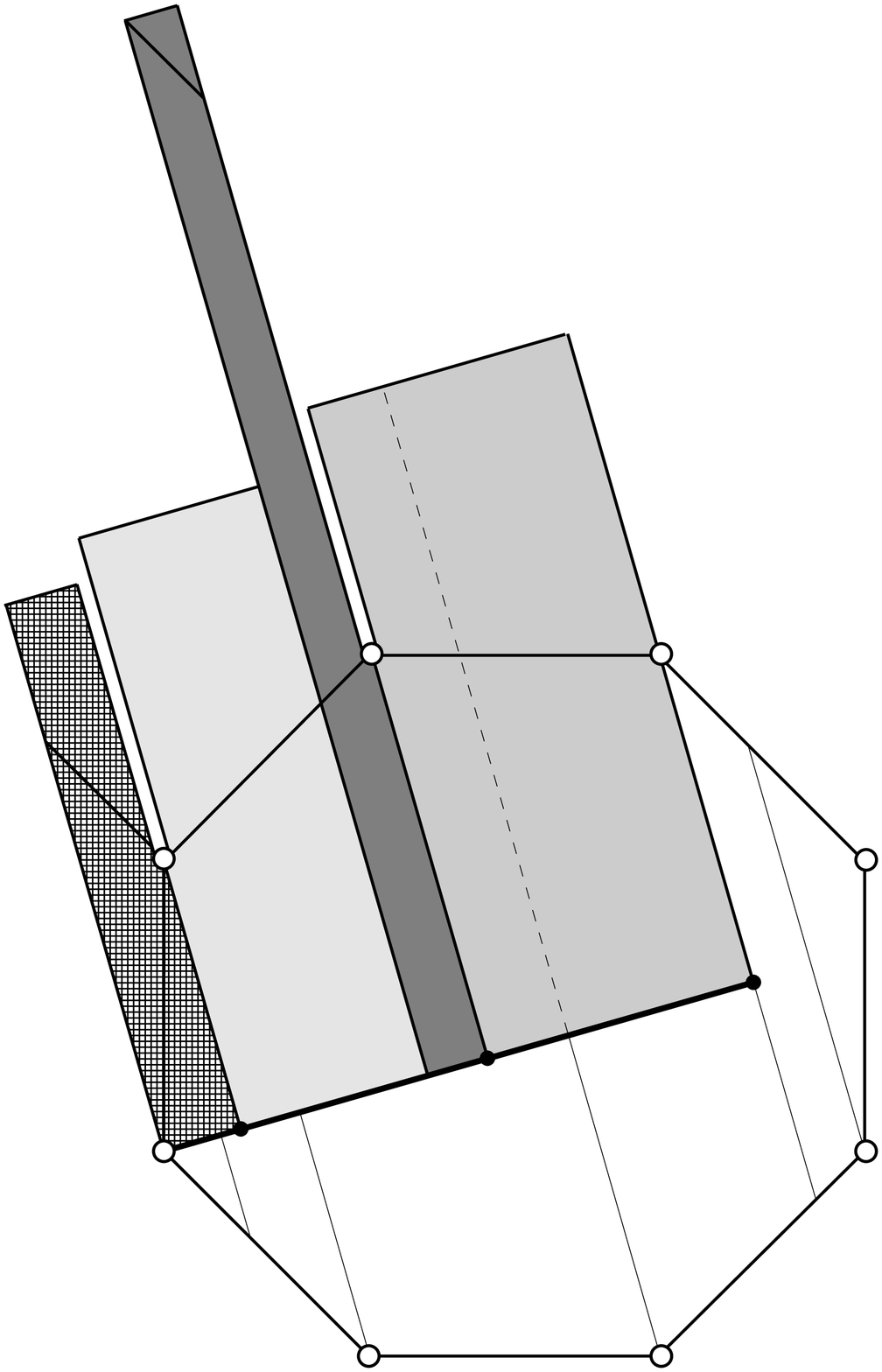}
\includegraphics{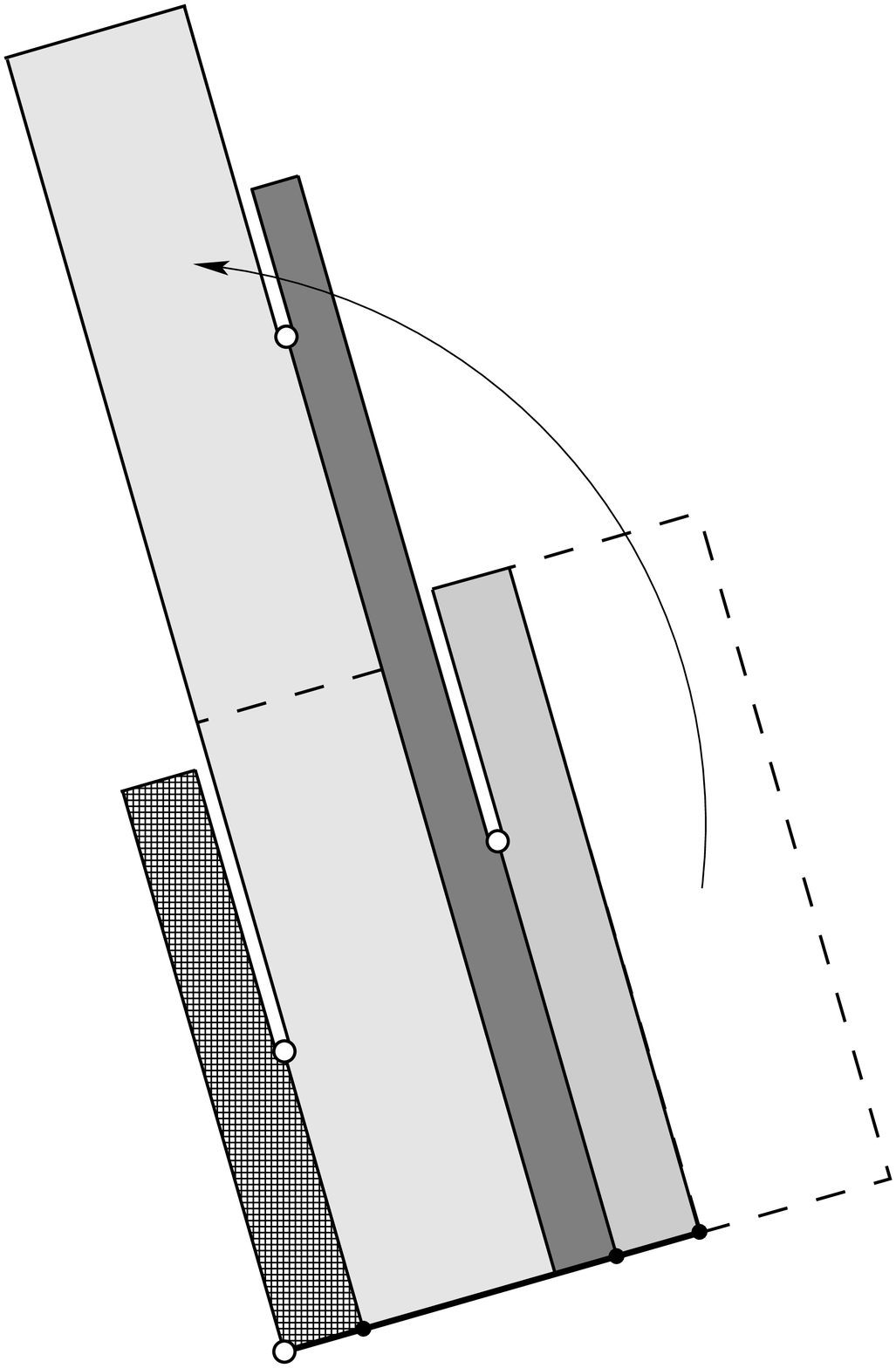}
\includegraphics{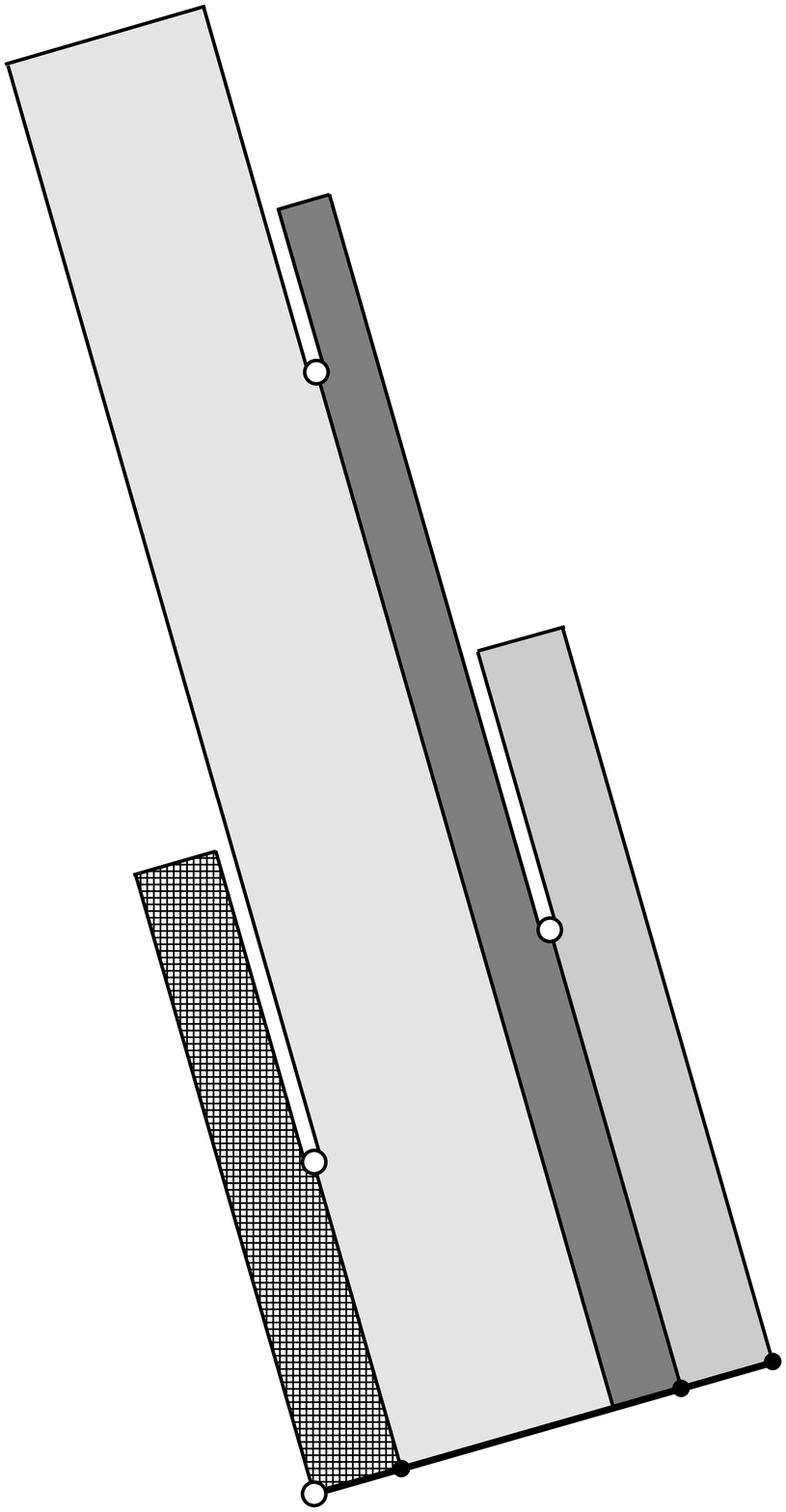}
\includegraphics{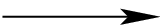}
\begin{picture}(100,50)(100,50) 
\put(230,-240){$\begin{pmatrix} e^{t_0} & 0 \\ 0 & e^{-t_0}\end{pmatrix}$}
\put(300,-70){$\underbrace{\hspace*{22truemm}}_{X'}$}
\put(360,-140){a)}
\put(160,-400){b)}
\put(345,-400){c)}
\put(111,-48){\Large $\vec v_1$}
\put(144,-9){\Large $\vec v_2$}
\put(196,-5){\Large $\vec v_3$}
\put(232,-37){\Large $\vec v_4$}
\put(116,-95){\Large $\vec v_4$}
\put(157,-126){\Large $\vec v_3$}
\put(204,-119){\Large $\vec v_2$}
\put(236,-84){\Large $\vec v_1$}
\end{picture}
\vspace{450bp}
\caption{
\label{zorich:fig:renormalization}
Idea  of  renormalization.  a)  Unwrap  the   flat  surface  into
``zippered rectangles''. b) Shorten the base of the corresponding
zippered rectangles.  c)Expand  the  resulting  tall  and  narrow
zippered rectangle  horizontally  and  contract  it vertically by
same factor $e^{t_0}$.
}
\end{figure}

The  reason  why the  deviation  of  the  cycles  $c_j$  from the
asymptotic direction  is  governed  by the Teichm\"uller geodesic
flow is illustrated  in Fig.~\ref{zorich:fig:renormalization}. In
a    sense,     we     follow     the     initial     ideas    of
H.~Masur~\cite{zorich:Masur:Annals:82}           and           of
W.~Veech~\cite{zorich:Veech:Annals:82}.

Fix a horizontal segment $X$ and emit a  vertical trajectory from
some point $x$ in $X$. When the trajectory intersects $X$ for the
first time join the corresponding  point  $T(x)$  to the original
point  $x$  along $X$ to obtain  a  closed loop. Here $T:X\to  X$
denotes the  first return map to  the transversal $X$  induced by
the vertical flow. Denote by  $c(x)$  the  corresponding cycle in
$H_1(S;\Z{})$. Let the interval  exchange  transformation $T:X\to
X$ decompose $X$ into $\noi$ subintervals $X_1\sqcup \dots \sqcup
X_\noi$. It is easy to see that the ``first return cycle'' $c(x)$
is piecewise  constant: we have $c(x)=c(x')=:c(X_j)$ whenever $x$
and   $x'$   belong  to   the   same   subinterval   $X_j$,   see
Fig.~\ref{zorich:fig:renormalization}. It is easy to see that
$$
c_N(x)=c(x)+c(T(x))+\dots+c(T^{N-1}(x))
$$

The average of this sum with respect to the ``time'' $N$ tends to
the asymptotic cycle $c$. We need to study the deviation  of this
sum  from  the value $N\cdot c$.  To  do this consider a  shorter
subinterval $X'$ as in Fig.~\ref{zorich:fig:renormalization}.
Its length  is chosen  in such way, that the  first return map of
the  vertical   flow   again   induces   an   interval   exchange
transformation $T':X'\to  X'$  of  $\noi$ subintervals. New first
return cycles $c'(X'_k)$ to the  interval  $X'$  are expressed in
terms of  the initial first  return cycles $c(X_j)$ by the linear
relations below; the  lengths $|X'_k|$ of subintervals of the new
partition  $X'=X'_1\sqcup\dots\sqcup X'_\noz$  are  expressed  in
terms  of  the  lengths  $|X_j|$ of subintervals of  the  initial
partition by dual linear relations:
$$
c'(X'_k)=\sum_{j=1}^\noi \Cocycle_{j k}\cdot c(X_j) \qquad\qquad
|X_j|=\sum_{k=1}^\noi \Cocycle_{j k}\cdot |X'_k|\,,
$$
where a  nonnegative integer matrix $\Cocycle_{jk}$ is completely
determined  by  the  initial  interval  exchange   transformation
$T:X\to X$ and by the choice of $X'\subset X$.

To construct the cycle $c_N$ representing a long piece of leaf of
the  vertical  foliation  we  followed the trajectory  $x,  T(x),
\dots,   T^{N-1}(x)$   of   the    initial    interval   exchange
transformation $T:X\to  X$ and computed the corresponding ergodic
sum. Passing to a shorter  horizontal  interval  $X'\subset X$ we
can  follow  the trajectory $x, T'(x), \dots, (T')^{N'-1}(x)$  of
the new  interval exchange transformation $T':X'\to X'$ (provided
$x\in X'$). Since  the subinterval $X'$  is shorter than  $X$  we
cover the  initial piece of trajectory of the  vertical flow in a
smaller number $N'$ of steps. In other words, passing from $T$ to
$T'$  we  accelerate  the  time:  it  is  easy  to see  that  the
trajectory  $x,   T'(x),   \dots,   (T')^{N'-1}(x)$  follows  the
trajectory  $x,  T(x), \dots, T^{N-1}(x)$ but jumps over  several
iterations of $T$ at a time.

This approach would not be efficient if the new first  return map
$T':X'\to X'$ would be more complicated than the initial one. But
we know that passing from $T$  to $T'$ we stay within a family of
interval exchange transformations of some fixed  number $\noi$ of
subintervals, and, moreover, that the new ``first return cycles''
$c'(X'_k)$ and the lengths  $|X'_k|$  of the new subintervals are
expressed  in  terms  of  the  initial  ones  by   means  of  the
$\noi\times\noi$-matrix  $\Cocycle$,  which  depends only on  the
choice of $X'\subset X$ and which can be easily computed.

Our strategy can be now formulated as follows. One can  define an
explicit  algorithm  (generalizing   Euclidean  algorithm)  which
canonically  associates  to  an interval exchange  transformation
$T:X\to X$ some  specific subinterval $X'\subset X$ and, hence, a
new interval  exchange transformation $T':X'\to X'$. Similarly to
the  Euclidean   algorithm   our  algorithm  is  invariant  under
proportional  rescaling of  $X$  and $X'$, so,  when  we find  it
convenient, we  can always rescale  the length of the interval to
one. This algorithm  can  be considered as a  map  $\cT$ from the
space of all  interval exchange transformations of a given number
$\noi$ of  subintervals  to  itself.  Applying  recursively  this
algorithm    we    construct    a   sequence   of    subintervals
$X=X^{(0)}\supset X^{(1)} \supset  X^{(2)}  \supset \dots $ and a
sequence      of      matrices       $\Cocycle=\Cocycle(T^{(0)}),
\Cocycle(T^{(1)}), \dots $ describing  transitions  form interval
exchange transformation $T^{(r)}:X^{(r)} \to X^{(r)}$ to interval
exchange  transformation   $T^{(r+1)}:X^{(r+1)}  \to  X^{(r+1)}$.
Taking    a    product     $\Cocycle^{(s)}=\Cocycle(T^{(0)})\cdot
\Cocycle(T^{(1)})\cdot  \dots\cdot  \Cocycle(T^{(s-1)})$  we  can
immediately express the ``first return cycles''  to a microscopic
subinterval $X^{(s)}$  in  terms  of  the  initial ``first return
cycles'' to $X$. Considering  now  the matrices $\Cocycle$ as the
values  of  a  matrix-valued  function on the space  of  interval
exchange transformations, we  realize  that we study the products
of matrices $\Cocycle$ along the orbits $T^{(0)}, T^{(1)}, \dots,
T^{(s-1)}$  of  the   map  on  the  space  of  interval  exchange
transformations. When the map is ergodic with respect to a finite
measure, the properties of these  products  are  described by the
Oseledets theorem, and the cycles $c_N$ have a deviation spectrum
governed by the Lyapunov  exponents  of the cocycle $\Cocycle$ on
the space of interval exchange transformations.

Note that the first return cycle  to  the  subinterval  $X^{(s)}$
(which is very short) represents the cycle $c_N$ corresponding to
a very long trajectory $x, T(x), ..., T^{N-1}(x)$  of the initial
interval   exchange   transformation.   In   other   words,   our
renormalization   procedure   $\cT$  plays  a  role  of  a   time
acceleration machine: morally, instead of getting the cycle $c_N$
by  following  a  trajectory  $x, T(x), ..., T^{N-1}(x)$  of  the
initial interval exchange transformation for the exponential time
$N\sim\exp(const\cdot s)$ we obtain the cycle $c_N$ applying only
$s$  steps  of the  renormalization  map $\cT$  on  the space  of
interval exchange transformations.

It  remains  to establish the relation between  the  the  cocycle
$\Cocycle$ over  the  map  $\cT$  and  the Teichm\"uller geodesic
flow. Conceptually, this relation was elaborated  in the original
paper of W.~Veech~\cite{zorich:Veech:Annals:82}.

First let us discuss how can one ``almost  canonically'' (that is
up   to   a  finite  ambiguity)  choose  a  zippered   rectangles
representation     of     a     flat    surface.    Note     that
Fig.~\ref{zorich:fig:renormalization}  suggests   the  way  which
allows   to   obtain   infinitely   many   zippered    rectangles
representations of the  same flat surface: we chop an appropriate
rectangle on the right,  put  it atop the corresponding rectangle
and then repeat the procedure  recursively.  This  resembles  the
situation  with  a  representation  of  a   flat   torus   by   a
parallelogram:   a   point   of   the   fundamental   domain   in
Fig.~\ref{zorich:fig:space:of:flat:tori}  provides  a   canonical
representative though any point of the corresponding $\SLZ$-orbit
represents  the   same   flat  torus.  A  ``canonical''  zippered
rectangles decomposition of a  flat  surface also belongs to some
fundamental  domain.   Following  W.~Veech  one  can  define  the
fundamental  domain  in  terms  of  some  specific  choice  of  a
``canonical'' horizontal  interval  $X$.  Namely, let us position
the left endpoint of $X$ at a conical singularity. Let  us choose
the  length  of  $X$  in  such  way that  the  interval  exchange
transformation  $T:X\to  X$ induced by the first  return  of  the
vertical flow to $X$ has minimal possible number $\noi=2g+\noz-1$
of  subintervals   under  exchange.  Among  all  such  horizontal
segments $X$ choose the shortest one,  which  length  is  greater
than or  equal to one.  This construction is applicable to almost
all flat surfaces; the finite ambiguity corresponds to the finite
freedom  in the  choice  of the conical  singularity  and in  the
choice of the horizontal ray adjacent to it.

Since the interval $X$ defines  a  decomposition  of (almost any)
flat     surface    into     ``zippered     rectangles''     (see
Fig.~\ref{zorich:fig:renormalization}) we can pass from the space
of flat  surfaces to the  space of zippered rectangles (which can
be considered as  a finite ramified  covering over the  space  of
flat surfaces).  Teichm\"uller  geodesic  flow lifts naturally to
the space of  zippered rectangles. It acts on zippered rectangles
by expansion in horizontal direction and  contraction in vertical
direction;  i.e.  the zippered rectangles are modified by  linear
transformations   $\begin{pmatrix}    e^t    &    0    \\   0   &
e^{-t}\end{pmatrix}$.  However,  as  soon  as  the  Teichm\"uller
geodesic  flow  brings   us  out  of  the
fundamental domain, we have to  modify  the  zippered  rectangles
decomposition to   the  ``canonical one'' corresponding to the
fundamental            domain.            (Compare             to
Fig.~\ref{zorich:fig:space:of:flat:tori} where the  Teichm\"uller
geodesic flow corresponds to the  standard  geodesic  flow in the
hyperbolic metric  on  the  upper  half-plane.) The corresponding
modification  of   zippered   rectangles   (chop  an  appropriate
rectangle on the  right, put it atop the corresponding rectangle;
repeat the procedure several times, if  necessary) is illustrated
in Fig.~\ref{zorich:fig:renormalization}.

Now everything is ready to establish  the  relation  between  the
Teichm\"uller geodesic flow  and the map  $\cT$ on the  space  of
interval exchange transformations.

Consider some codimension one subspace $\Upsilon$ in the space of
zippered  rectangles  transversal  to the Teichm\"uller  geodesic
flow. Say, $\Upsilon$ might be  defined  by  the requirement that
the base  $X$ of the  zippered rectangles decomposition has length one, $|X|=1$.
This    is    the    choice    in   the   original    paper    of
W.~Veech~\cite{zorich:Veech:Annals:82};    under    this   choice
$\Upsilon$ represents part of the  boundary  of  the  fundamental
domain  in   the  space  of  zippered  rectangles.  Teichm\"uller
geodesic     flow     defines    the     first     return     map
$\cS:\Upsilon\to\Upsilon$ to the section $\Upsilon$. The map $\cS$ can
be  described  as follows.  Take  a  flat  surface  of  unit area
decomposed  into  zippered rectangles $Z$ with the  base  $X$  of
length  one.   Apply   expansion   in  horizontal  direction  and
contraction in vertical direction. For some $t_0(Z)$ the deformed
zippered     rectangles     can     be    rearranged    as     in
Fig.~\ref{zorich:fig:renormalization}  to  get back to  the  base  of
length one;  the result is the image of  the map $\cS$. Actually,
we    can     first     apply    the    rearrangement    as    in
Fig.~\ref{zorich:fig:renormalization}  to  the  initial  zippered
rectangles $Z$ and then apply the transformation $\begin{pmatrix}
e^{t_0} & 0 \\ 0 &  e^{-t_0}\end{pmatrix}$ ---  the  two  operations
commute.  This  gives,  in  particular, an explicit  formula  for
$t_0(Z)$. Namely let $|X_\noi|$ be the  width  of  the  rightmost
rectangle and  let $|X_k|$ be the  width of the  rectangle, which
top horizontal side  is glued to  the rightmost position  at  the
base $X$.  (For  the  upper  zippered  rectangle decomposition in
Fig.~\ref{zorich:fig:renormalization} we  have $n=4$ and  $k=2$.)
Then
$$
t_0=-\log\big(1-\min(|X_n|,|X_k|)\big).
$$

Recall  that  a  decomposition  of a flat surface  into  zippered
rectangles naturally defines an interval exchange  transformation
--- the first return map of  the vertical flow to the base $X$ of
zippered  rectangles.  Hence,  the  map  $\cS$  of  the  subspace
$\Upsilon$ of zippered rectangles  defines  an induced map on the
space  of  interval exchange transformations. It remains to  note
that this  induced map is exactly the map  $\cT$. In other words,
the map $\cS:\Upsilon\to\Upsilon$ induced by the  first return of
the Teichm\"uller geodesic flow to  the  subspace  $\Upsilon$  of
zippered rectangles is  the suspension of  the map $\cT$  on  the
space of interval exchange transformations.

We complete with a remark concerning the choice of a section. The
natural  section  $\Upsilon$ chosen  in  the  original  paper  of
W.~Veech~\cite{zorich:Veech:Annals:82} is in a  sense  too large:
the corresponding invariant measure (induced from  the measure on
the space of  flat surfaces) is infinite. Choosing an appropriate
subset $\Upsilon'\subset\Upsilon$  one  can  get finite invariant
measure. Moreover, the subset  $\Upsilon'$  can be chosen in such
way     that     the     corresponding    first    return     map
$\cS':\Upsilon'\to\Upsilon'$ of  the Teichm\"uller geodesic  flow
is  a suspension of  some  natural  map  $\cG$ on  the  space  of
interval                exchange                 transformations,
see~\cite{zorich:Zorich:Gauss:map}. According  to the results  of
H.~Masur~\cite{zorich:Masur:Annals:82}                        and
W.~Veech~\cite{zorich:Veech:Annals:82} the Teichm\"uller geodesic
flow is ergodic which implies ergodicity of  the maps $\cS'$
and  $\cG$.  To  apply  Oseledets theorem one  should,  actually,
consider the induced cocycle $B$ over this new  map $\cG$ instead
of the cocycle $\Cocycle$ over the map $\cT$ described above.

\section{Closed Geodesics on Flat Surfaces}
\label{zorich:s:Closed:Geodesics:on:Flat:Surfaces}

Consider  a flat  surface  $S$; we always  assume  that the  flat
metric on $S$ has trivial holonomy, and that the surface  $S$ has
finite  number  of cone-type singularities. By convention a  flat
surface is  endowed with a choice of direction,  refereed to as a
``vertical direction'', or as a ``direction to the North''. Since
the  flat  metric  has  trivial holonomy, this direction  can  be
transported in a unique way to any point of the surface.

A geodesic segment  joining  two conical singularities and having
no conical points in its  interior  is  called \textit{saddle connection}.
The case  when boundaries of  a saddle connection coincide is not
excluded:  a  saddle  connection  might join a conical  point  to
itself. In this section we  study  saddle  connections and closed
regular geodesics on a generic flat surface $S$ of genus $g\ge 2$.
In particular, we count them and
we explain the following curious phenomenon: saddle connections and
closed regular geodesics
often appear in pairs, triples, etc  of parallel saddle
connections  (correspondingly  closed  regular geodesics) of  the
same  direction  and length. When all saddle connections  (closed
regular  geodesics)   in   such   configuration   are  short  the
corresponding  flat  surface is almost degenerate; it is  located
close  to the  boundary  of the moduli  space.  A description  of
possible configurations  of  parallel  saddle connections (closed
geodesics)  gives   us  a  description  of  the  multidimensional
``cusps'' of the strata.

The results of this section are based on the joint work with
A.~Eskin and H.~Masur~\cite{zorich:Eskin:Masur:Zorich} and on
their work~\cite{zorich:Masur:Zorich}. A series of beautiful results
developing the counting problems considered here were recently
obtained by Ya.~Vorobets~\cite{zorich:Vorobets:uniform:bounds}.

\paragraph{Counting closed geodesics and saddle connections.}

Closed  geodesics  on flat surfaces of higher  genera  have  some
similarities  with ones  on  the torus. Suppose  that  we have  a
regular  closed  geodesic  passing  through a point  $x_0\in  S$.
Emitting a geodesic from a nearby point $x$ in the same direction
we obtain a parallel closed  geodesic  of the same length as  the
initial  one.  Thus,  closed  geodesics  appear  in  families  of
parallel closed geodesics.  However, in the torus case every such
family fills the entire torus while a family  of parallel regular
closed geodesics  on a flat surfaces  of higher genus  fills only
part of the surface. Namely,  it  fills a flat cylinder having  a
conical singularity on each of its boundaries.
Typically, a maximal cylinder of
closed regular geodesics is bounded by a pair of closed saddle
connections. Reciprocally, any saddle connection joining a conical
point $P$ to itself and coming back to $P$ at the angle $\pi$
bounds a cylinder filled with closed regular geodesics.

A geodesic representative  of  a homotopy class of  a  curve on a
flat surface is realized in general by a broken line  of geodesic
segments  with  vertices at  conical  points.  By  convention  we
consider  only  closed  \textit{regular}   geodesics   (which  by
definition  do  not  pass  through  conical   points)  or  saddle
connections (which by definition do  not  have  conical points in
its interior). Everywhere in  this  section we normalize the area
of a flat surface to one.

Let $N_{sc}(S,L)$ be the  number  of saddle connections of length
at most $L$  on  a  flat surface $S$. Let  $N_{cg}(S,L)$  be  the
number of maximal cylinders filled with  closed regular geodesics
of length at most $L$ on  $S$. It was proved by H.~Masur that for
any  flat  surface $S$  both  counting  functions  $N(S,L)$  grow
quadratically   in   $L$.   Namely,    there    exist   constants
$0<const_1(S)<const_2(S)<\infty$ such that
$$
const_1(S) \le N(S,L)/L^2 \le const_2(S)
$$
for $L$  sufficiently  large.  Recently Ya.~Vorobets has obtained
uniform estimates for the constants $const_1(S)$ and $const_2(S)$
which    depend     only     on     the     genus     of     $S$,
see~\cite{zorich:Vorobets:uniform:bounds}.      Passing      from
\emph{all} flat surfaces to \emph{almost all} surfaces in a given
connected  component of  a  given stratum one  gets  a much  more
precise result, see~\cite{zorich:Eskin:Masur}:

\begin{NNTheorem}[A.~Eskin and H.~Masur]
For almost  all flat surfaces  $S$ in any
stratum  $\cH(d_1,\dots,d_\noz)$  the  counting  functions
$N_{sc}(S,L)$ and $N_{cg}(S,L)$ have exact quadratic asymptotics
$$
\lim_{L\to\infty}\cfrac{N_{sc}(S,L)}{\pi L^2}=c_{sc}(S) \qquad
\lim_{L\to\infty}\cfrac{N_{cg}(S,L)}{\pi L^2}=c_{cg}(S)\
$$
Moreover, the Siegel--Veech constants $c_{sc}(S)$ (correspondingly $c_{cg}(S)$)
coincide
for   almost    all    flat    surfaces $S$   in
each connected   component
$\cH^{comp}_1(d_1,\dots,d_\noz)$ of the stratum.
\end{NNTheorem}

\paragraph{Phenomenon of higher multiplicities.}

Note that the direction to  the  North is well-defined even at  a
conical point of  a  flat surface, moreover, at  a  conical point $P_1$
with a cone angle $2\pi  k$  we have $k$ different directions  to
the North! Consider some saddle  connection  $\gamma_1=[P_1  P_2]$ with an
endpoint at $P_1$.  Memorize its direction,  say, let it  be  the
North-West  direction.  Let us launch a geodesic  from  the  same
starting  point  $P_1$ in one of the  remaining  remaining  $k-1$
North-West directions. Let us study how big is the chance  to hit
$P_2$ ones again,  and how  big is  the  chance to  hit it  after
passing the  same distance as before. We do  not exclude the case
$P_1=P_2$. Intuitively  it is clear  that the answer to the first
question is: ``the  chances  are low'' and to  the  second one is
``the chances are  even lower''. This makes the following theorem
(see~\cite{zorich:Eskin:Masur:Zorich}) somehow counterintuitive:

\begin{Theorem}[A.~Eskin, H.~Masur, A.~Zorich]
For almost any flat  surface $S$ in any stratum and for  any pair
$P_1,  P_2$  of   conical   singularities  on  $S$  the  function
$N_2(S,L)$  counting  the number  of  pairs  of  parallel  saddle
connections of  the same length  joining $P_1$ to $P_2$ has exact
quadratic asymptotics
$$
\lim_{L\to\infty}\cfrac{N_2(S,L)}{\pi L^2}= c_2>0\ ,
$$
where the Siegel--Veech constant $c_2$ depends only  on the connected
component  of  the stratum and on  the  cone angles at $P_1$  and
$P_2$.

For almost all flat surfaces $S$  in any stratum one cannot  find
neither a  single pair of parallel  saddle connections on  $S$ of
different   length,   nor  a  single  pair  of  parallel   saddle
connections joining different pairs of singularities.
\end{Theorem}

Analogous statements
(with some reservations for specific connected components of certain strata)
can be  formulated  for arrangements of $3, 4, \dots$
parallel saddle connections. The situation with
closed  regular  geodesics  is  similar: they might  appear
(also with some exceptions for specific connected components of certain strata)
in
families of  $2, 3, \dots, g-1$  distinct maximal cylinders
filled with parallel  closed  regular geodesics of equal length.
A general formula  for  the  Siegel--Veech constant in
the corresponding quadratic asymptotics is
presented at  the end of this section, while here
we  want  to  discuss  the numerical values  of  Siegel--Veech
constants in a simple concrete example. We consider the   principal   strata
$\cH(1,\dots,1)$ in small genera. Let $N_{k\underline\ cyl}(S,L)$
be the corresponding counting  function,  where $k$ is the number
of distinct maximal cylinders filled with parallel closed regular
geodesics of equal length bounded by $L$. Let

$$
c_{k\underline\ cyl}=\lim_{L\to\infty}\cfrac{N_{k\underline\ cyl}(S,L)}{\pi L^2}
$$
The table below (extracted from~\cite{zorich:Eskin:Masur:Zorich})
presents the  values of $c_{k\underline\ cyl}$ for $g=1,\dots,4$.
Note  that  for  a  generic  flat  surface  $S$  of genus  $g$  a
configuration of  $k\ge g$ cylinders  is not realizable, so we do
not fill the corresponding entry.

$$
\begin{array}{|c|c|c|c|c|}
\hline &&&&\\
[-\halfbls]k & g=1 & g=2 & g=3 & g=4 \\
\hline &&&&\\
[-\halfbls]&&&&\\
[-2\halfbls]\ 1\ &
\cfrac{1}{2}\cdot\cfrac{1}{\zeta(2)}\approx 0.304 &
\cfrac{5}{2}\cdot\cfrac{1}{\zeta(2)}\approx 1.52 &
\cfrac{36}{7}\cdot\cfrac{1}{\zeta(2)}\approx 3.13 &
\cfrac{3150}{377}\cdot\cfrac{1}{\zeta(2)}\approx 5.08 \\
[-1\halfbls]&&&&\\
\hline &&&&\\
[-\halfbls]&&&&\\
[-2\halfbls]\ 2\ &
- & - &
\cfrac{3}{14}\cdot\cfrac{1}{\zeta(2)}\approx 0.13 &
\cfrac{90}{377}\cdot\cfrac{1}{\zeta(2)}\approx 0.145 \\
[-\halfbls]&&&&\\
\hline &&&&\\
[-\halfbls]&&&&\\
[-2\halfbls]\ 3\ &
- & - & - &
\cfrac{5}{754}\cdot\cfrac{1}{\zeta(2)}\approx 0.00403 \\
[-\halfbls]&&&&\\
\hline
\end{array}
$$

Comparing these  values we see,  that our intuition was not quite
misleading. Morally, in genus $g=4$  a  closed  regular  geodesic
belongs to  a  one-cylinder family with ``probability'' $97.1\%$,
to a two-cylinder family with  ``probability''  $2.8\%$  and to a
three-cylinder family  with  ``probability''  only $0.1\%$ (where
``probabilities'' are calculated  proportionally to the Siegel--Veech
constants $5.08\,:\,0.145\,:\,0.00403$).

\paragraph{Rigid configurations of saddle connections
and ``cusps'' of the strata.}

A  saddle  connection or  a  regular closed  geodesic  on a  flat
surface $S$ persists under  small  deformations of $S$ inside the
corresponding stratum. It  might happen that any deformation of a
given flat surface which shortens some specific saddle connection
necessarily shortens some other saddle connections. We say that a
collection $\{\gamma_1, \dots, \gamma_n\}$ of saddle  connections
is  \textit{rigid}  if any sufficiently small deformation of  the
flat  surface   inside  the  stratum  preserves  the  proportions
$|\gamma_1|:|\gamma_2|: \dots:|\gamma_n|$ of the  lengths  of all
saddle   connections    in   the   collection.   It   was   shown
in~\cite{zorich:Eskin:Masur:Zorich} that all  saddle  connections
in any  rigid  collection  are  \textit{homologous}.  Since their
directions and lengths can be expressed in terms  of integrals of
the holomorphic  1-form  $\omega$ along corresponding paths, this
implies  that  homologous saddle  connections  $\gamma_1,  \dots,
\gamma_n$ are  parallel and have equal  length and either  all of
them  join the  same  pair of distinct  singular  points, or  all
$\gamma_i$ are closed loops.

\begin{figure}
%
 %
 %
 \includegraphics{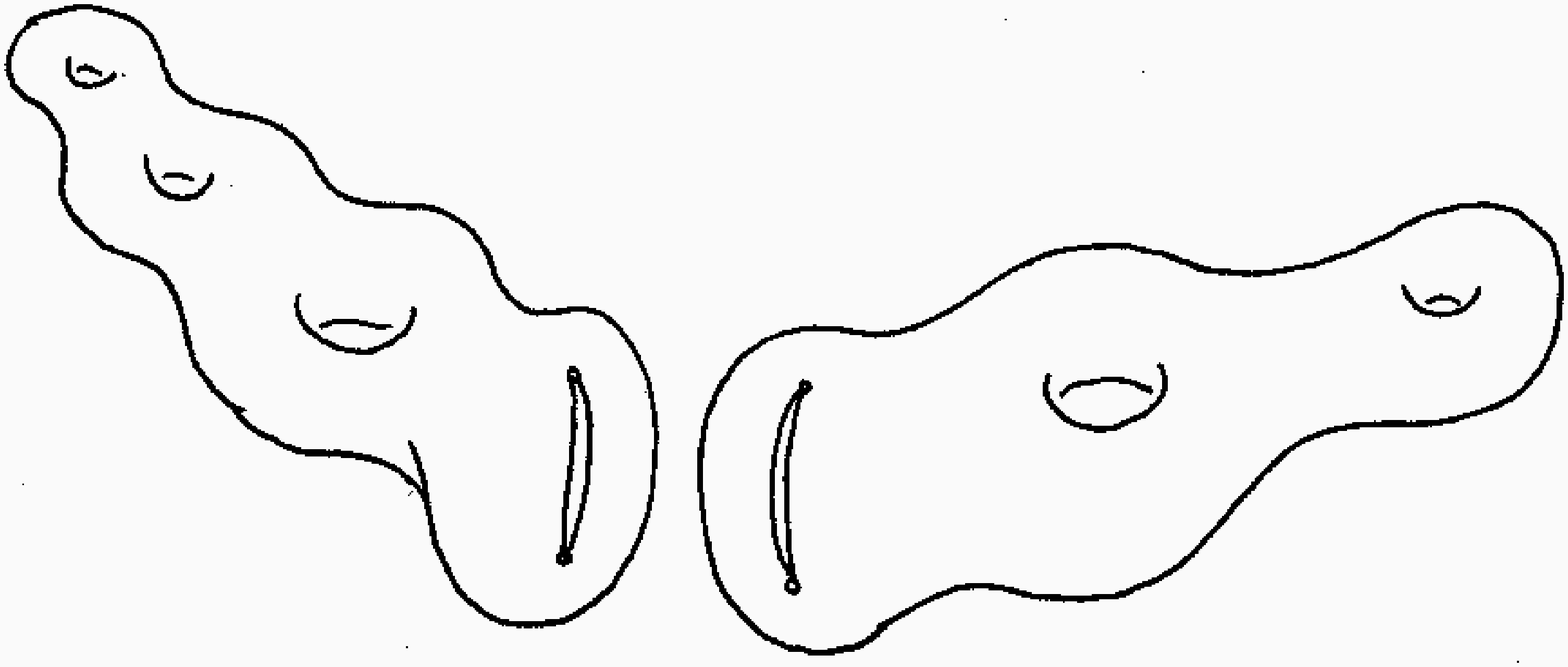}
 \includegraphics{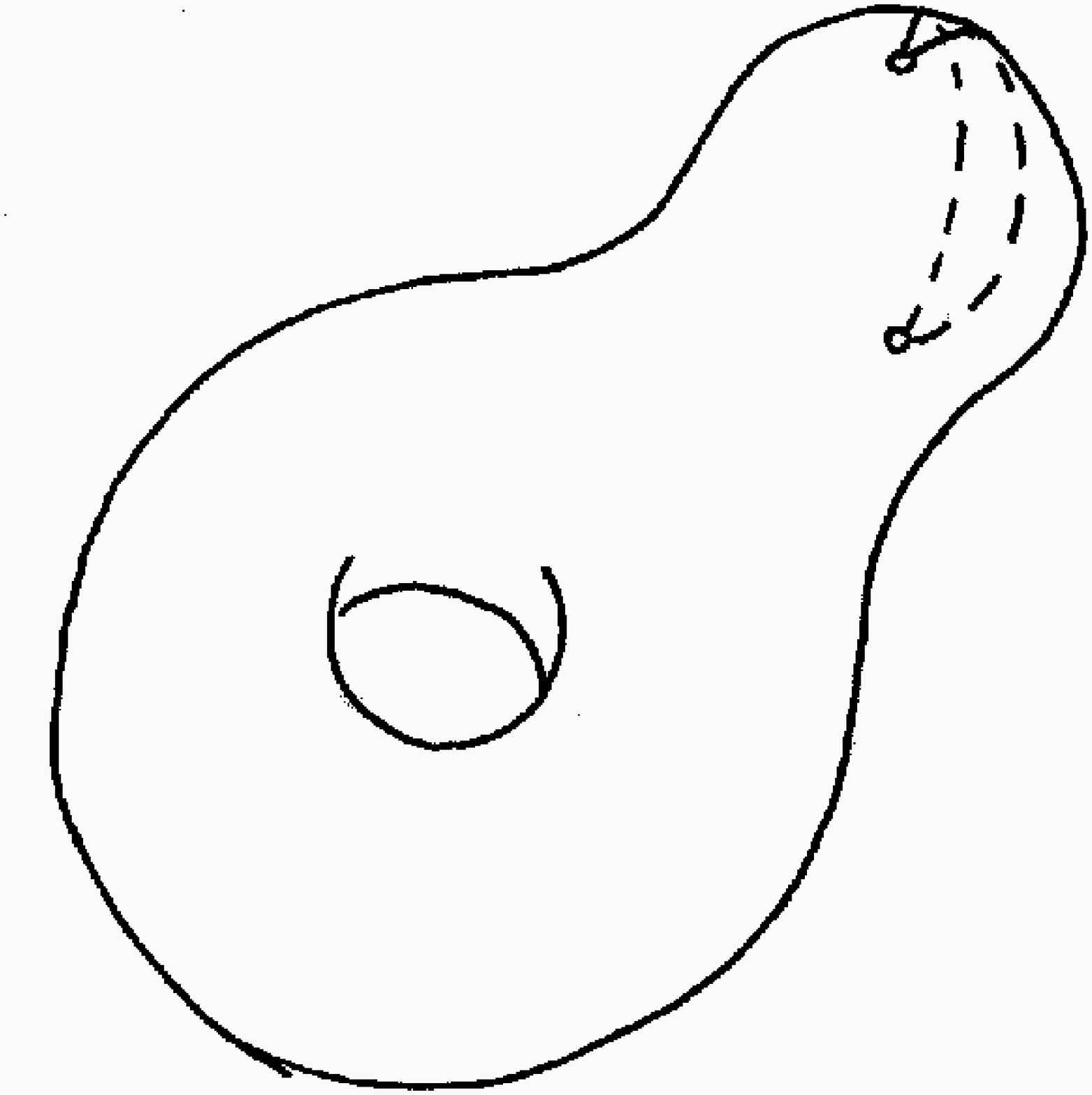}
 %
 %
 \includegraphics{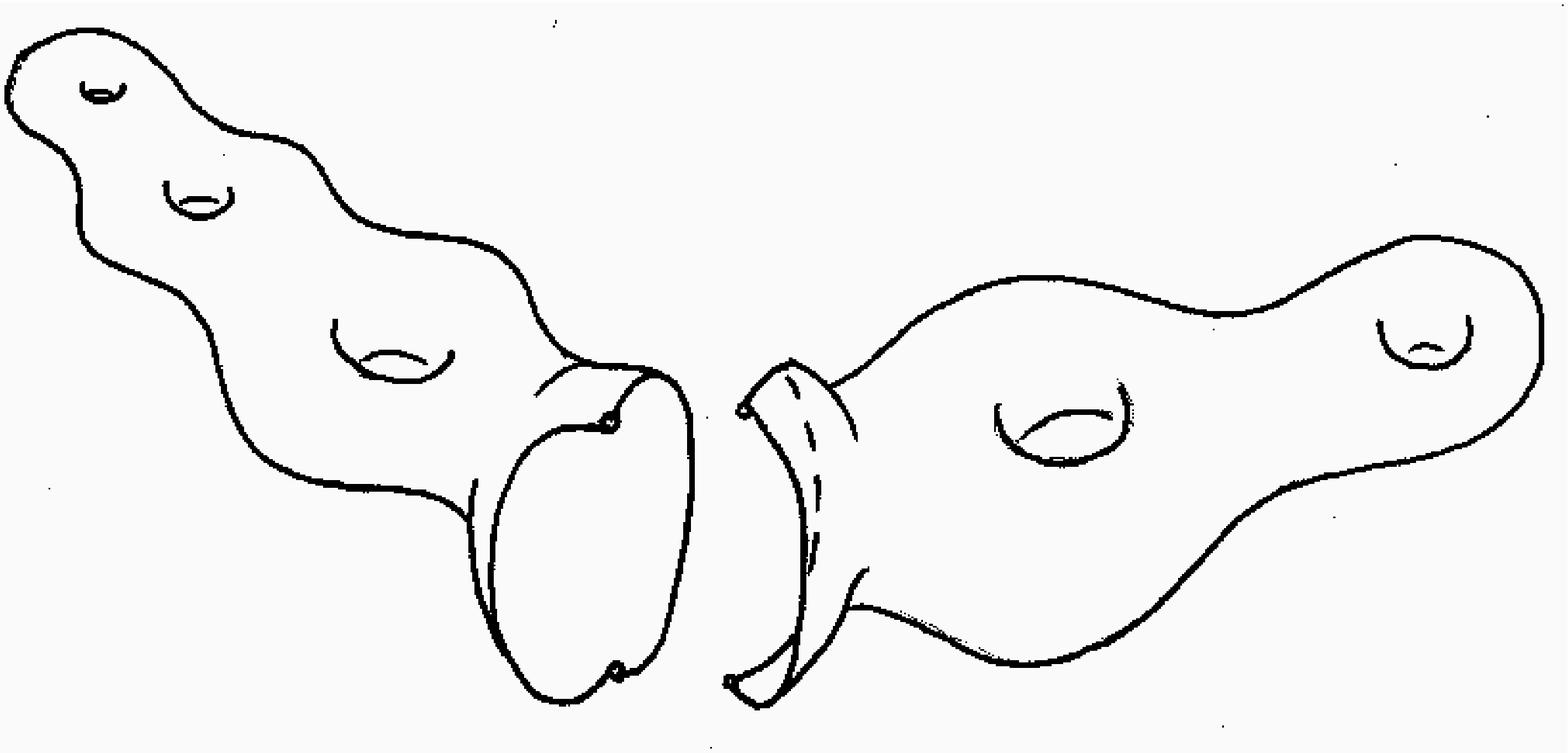}
 \includegraphics{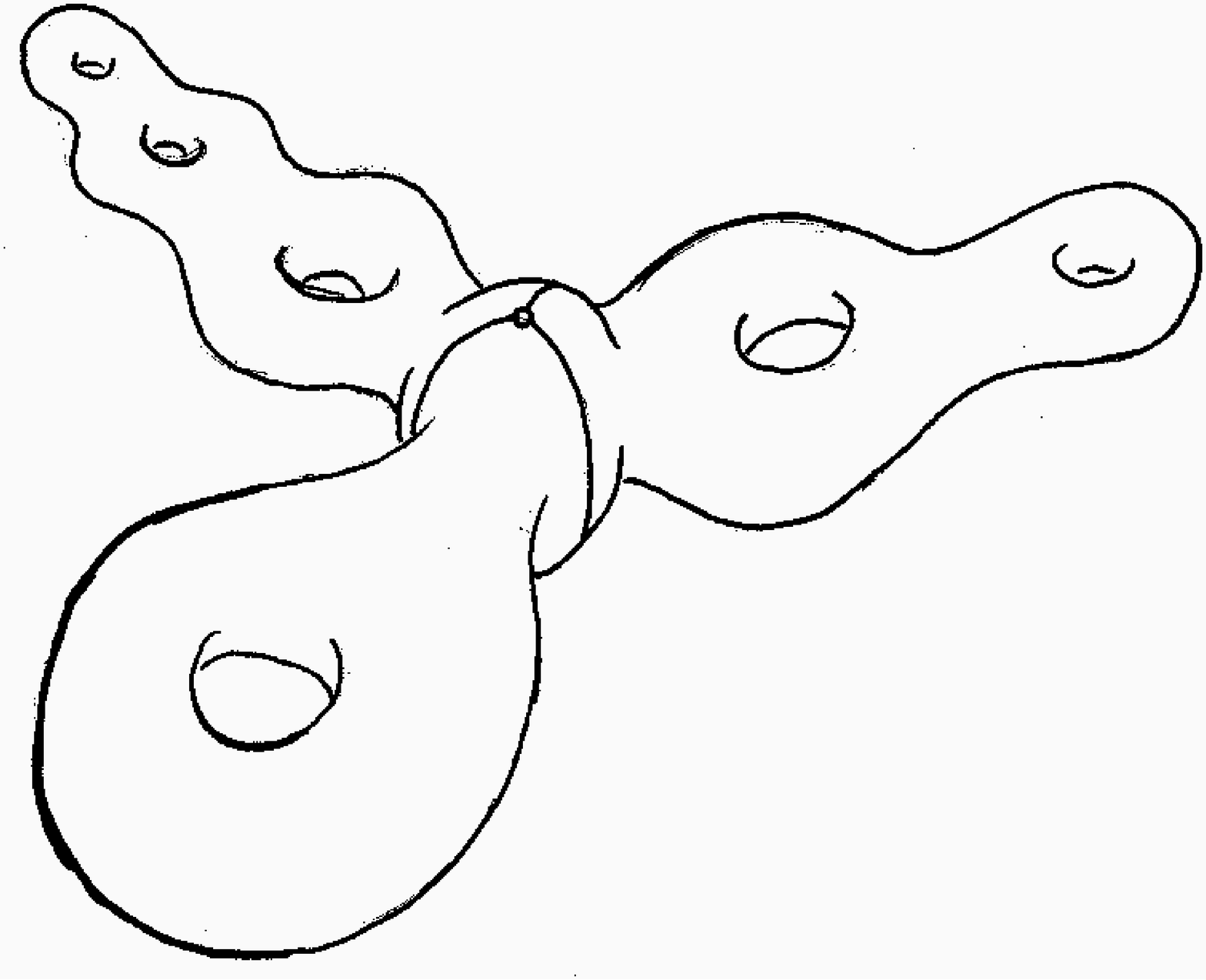}
%
%
\begin{picture}(-160,0)(-160,0)
\put(0,0){
 \begin{picture}(0,0)(0,0) 
 \put(-153,-36){$S'_3$}
 \put(-6,-46){$S'_2$} 
 \end{picture}
}
\put(5,0){
 \begin{picture}(0,0)(0,0) 
 \put(-94,-127){$S'_1$}
 \end{picture}
}
\put(0,0){
 \begin{picture}(0,0)(0,0) 
 \put(20,-12){$S_3$}
 \put(155,-22){$S_2$}
 \end{picture}
}
\put(0,0){
 \begin{picture}(0,0)(0,0) 
 \put(80,-180){$S_1$}
 \put(20,-117){$S_3$}
 \put(145,-126){$S_2$}
 \put(84,-140){$\scriptstyle \gamma_1$}
 \put(70,-140){$\scriptstyle \gamma_2$}
 \put(79,-125){$\scriptstyle \gamma_3$}
 \end{picture}
}
\end{picture}
\vspace{195bp}
\caption{
\label{zorich:fig:dist:sad:mult}
Multiple homologous saddle connections, topological picture
(after~\cite{zorich:Eskin:Masur:Zorich})}
\end{figure}

This implies that when  saddle  connections in a rigid collection
are contracted by a continuous  deformation,  the  limiting  flat
surface generically decomposes into several connected  components
represented by  nondegenerate  flat surfaces $S'_1, \dots, S'_k$,
see Fig.~\ref{zorich:fig:dist:sad:mult},
where  $k$ might vary  from  one  to  the genus  of  the  initial
surface.  Let  the  initial  surface  $S$  belong  to  a  stratum
$\cH(d_1,\dots,d_\noz)$.  Denote  the   set  with  multiplicities
$\{d_1,\dots,d_\noz\}$  by  $\beta$.  Let $\cH(\beta'_j)$ be  the
stratum      ambient      for       $S'_j$.      The      stratum
$\cH(\beta')=\cH(\beta'_1)\sqcup \dots \sqcup  \cH(\beta'_k)$  of
disconnected  flat  surfaces  $S'_1\sqcup  \dots\sqcup  S'_k$  is
referred  to  as  a  \textit{principal boundary} stratum  of  the
stratum $\cH(\beta)$. For any connected component  of any stratum
$\cH(\beta)$ the paper~\cite{zorich:Eskin:Masur:Zorich} describes
all  principal  boundary  strata;  their  union   is  called  the
principal  boundary   of  the  corresponding   connected
component of $\cH(\beta)$.

The  paper~\cite{zorich:Eskin:Masur:Zorich}   also  presents  the
inverse construction. Consider any flat surface $S'_1\sqcup \dots
\sqcup  S'_k\in  \cH(\beta')$ in the principal boundary  of
$\cH(\beta)$;  consider  a sufficiently small value of a  complex
parameter  $\varepsilon\in\C{}$.  One  can  reconstruct the  flat
surface $S\in\cH(\beta)$  endowed with a collection of homologous
saddle  connections   $\gamma_1,   \dots,   \gamma_n$  such  that
$\int_{\gamma_i}\omega=\varepsilon$, and such  that  degeneration
of  $S$  contracting  the  saddle connections $\gamma_i$  in  the
collection gives the surface $S'_1\sqcup \dots \sqcup S'_k$. This
inverse construction  involves  several  surgeries  of  the  flat
structure. Having  a  disconnected flat surface $S'_1\sqcup \dots
\sqcup S'_k$ one applies an  appropriate  surgery  to each $S'_j$
producing a surface $S_j$  with  boundary. The surgery depends on
the  parameter  $\varepsilon$: the  boundary  of  each  $S_j$  is
composed from  two  geodesic segments of lengths $|\varepsilon|$;
moreover, the  boundary  components  of  $S_j$  and $S_{j+1}$ are
compatible, which allows to glue  the  compound  surface $S$ from
the    collection    of    surfaces     with     boundary,    see
Fig.~\ref{zorich:fig:dist:sad:mult} as an example.

A collection $\gamma=\{\gamma_1, \dots, \gamma_n\}$ of homologous
saddle connections determines the following data on combinatorial
geometry of the decomposition $S\setminus \gamma$:  the number of
components, their  boundary  structure,  the singularity data for
each  component,  the cyclic order in which  the  components  are
glued   to   each   other.   These  data  are   referred   to  as
a \emph{configuration} of  homologous   saddle   connections.   A
configuration   $\cC$   uniquely  determines   the  corresponding
boundary stratum $\cH(\beta'_{\cC})$;  it  does not depend on the
collection $\gamma$ of homologous saddle connections representing
the configuration $\cC$.

The constructions  above explain how configurations $\cC$ of homologous
saddle connections  on  flat  surfaces $S\in\cH(\beta)$ determine
the ``cusps''  of  the  stratum  $\cH(\beta)$.  Consider a subset
$\cH_1^{\varepsilon}(\beta)\subset\cH(\beta)$ of surfaces of area
one having a  saddle connection shorter than $\varepsilon$. Up to
a  subset  $\cH_1^{\varepsilon,thin}(\beta)$ of  negligibly small
measure
the    set     $\cH_1^{\varepsilon,thick}(\beta)     =
\cH_1^{\varepsilon}(\beta)                              \setminus
\cH_1^{\varepsilon,thin}(\beta)$  might   be  represented  as   a
disjoint union
$$
\cH_1^{\varepsilon,thick}(\beta)\approx\bigsqcup_\cC\cH_1^\varepsilon(\cC)
$$
of    neighborhoods $\cH_1^\varepsilon(\cC)$   of    the   corresponding
``cusps'' $\cC$. Here $\cC$ runs over a finite set of
configurations admissible  for the given
stratum $\cH_1(\beta)$; this set is explicitly described in~\cite{zorich:Eskin:Masur:Zorich}.

When a configuration $\cC$ is  composed  from  homologous  saddle
connections  joining  \textit{distinct} zeroes,  the neighborhood
$\cH_1^\varepsilon(\cC)$ of the cusp $\cC$ has the structure of a
fiber   bundle   over   the   corresponding   boundary    stratum
$\cH(\beta'_{\cC})$ (up  to  a
difference in a set of a negligibly small measure). A fiber of
this bundle is represented by a finite
cover over the Euclidean  disc of radius $\varepsilon$
ramified at the center of the disc.
Moreover, the
canonical  measure  in  $\cH_1^\varepsilon(\cC)$ decomposes into  a
product measure of  the canonical measure in the boundary stratum
$\cH(\beta'_{\cC})$  and  the  Euclidean  measure  in  the  fiber
(see~\cite{zorich:Eskin:Masur:Zorich}), so
\begin{equation}
\label{zorich:eq:volume:of:a:cusp}
\Vol\left(\cH_1^\varepsilon(\cC)\right)=
(\text{combinatorial factor})\cdot\pi\varepsilon^2\cdot
\prod_{j=1}^k\Vol\cH_1(\beta'_j)\ +\ o(\varepsilon^2) .
\end{equation}

\begin{NNRemark}
We   warn    the   reader   that   the   correspondence   between
compactification of the moduli space of Abelian differentials and
the Deligne---Mamford compactification of  the  underlying moduli
space  of  curves is  not  straightforward.  In  particular,  the
desingularized stable  curve  corresponding  to the limiting flat
surface generically \textit{is  not}  represented as a union of
Riemann surfaces corresponding to $S'_1, \dots, S'_k$ --- the stable
curve might contain more components.
\end{NNRemark}

\paragraph{Evaluation of the Siegel--Veech constants.}

Consider a flat surface $S$. To  every  closed  regular  geodesic
$\gamma$ on $S$  we can associate  a vector $\vec  v(\gamma)$  in
$\R{2}$ having the length and the direction of $\gamma$. In other
words, $\vec v  = \int_\gamma\omega$, where we consider a complex
number  as   a   vector   in   $\R{2}\simeq\C{}$.  Applying  this
construction to all closed regular geodesic on $S$ we construct a
discrete set $V(S)\subset \R{2}$. Consider the following operator
$  f  \mapsto  \hat{f}$  from functions with compact  support  on
$\R{2}$    to    functions    on     a     connected    component
$\cH_1^{comp}(\beta)$           of           the          stratum
$\cH_1(\beta)=\cH_1(d_1,\dots,d_\noz)$:
$$
\hat{f}(S):=\sum_{\vec v\in V(S)} f(\vec v)
$$
Function $\hat{f}(S)$ generalizes the counting function $N_{cg}(S,L)$
introduced in the beginning of this section. Namely, when $f=\chi_L$
is the characteristic function $\chi_L$ of the disc of
radius  $L$ with the center at the origin of $\R{2}$,
the function $\hat\chi_L(S)$ counts the number of regular closed geodesics
of length at most $L$ on a flat surface $S$.

\begin{NNTheorem}[W.~Veech]
For any function $f:\R{2}\to\R{}$ with compact support the
following equality is valid:
\begin{equation}
\label{zorich:eq:Veech:counting:theorem}
\cfrac{1}{\Vol\cH^{comp}_1(\beta)}
\int_{\cH^{comp}_1(\beta)} \hat{f}(S)\, d\nu_1 =
C \int_{\R{2}} f(x,y)\, dx\, dy\ ,
\end{equation}
where the constant $C$ does not depend on the function $f$.
\end{NNTheorem}

Note  that  this is  an  exact  equality.  In  particular,
choosing the characteristic function $\chi_L$ of a disc of
radius  $L$ as a function $f$
we  see  that   for   any   positive  $L$  the
\emph{average} number of closed regular geodesics not longer than
$L$ on flat surfaces $S\in\cH^{comp}_1(\beta)$ is  exactly
$C\cdot\pi  L^2$,  where the Siegel--Veech constant $C$ does  not
depend   on   $L$,  but   only   on   the   connected   component
$\cH^{comp}_1(\beta)$.

The theorem  of  Eskin  and Masur cited
above   tells   that   for   large   values   of  $L$  one   gets
approximate equality
$\hat{\chi}_L(S)\approx c_{cg}\cdot \pi L^2$ ``pointwisely''  for
almost        all        individual         flat         surfaces
$S\in\cH^{comp}_1(d_1,\dots,d_\noz)$.      It      is      proved
in~\cite{zorich:Eskin:Masur} that the corresponding Siegel--Veech
constant  $c_{cg}$   coincides   with   the   constant   $C$   in
equation~\eqref{zorich:eq:Veech:counting:theorem} above.

Actually,  the  same  technique  can  be  applied to count
separately  pairs,  triples,  or  any other specific configurations
$\cC$ of homologous saddle connections. Every time when we find a
collection  of  homologous  saddle connections $\gamma_1,  \dots,
\gamma_n$  representing  the chosen
configuration    $\cC$    we    construct    a    vector    $\vec
v=\int_{\gamma_i}\omega$. Since all $\gamma_1,\dots,\gamma_n$ are
homologous, we can  take any of  them as $\gamma_i$.  Taking  all
possible collections  of  homologous  saddle  connections  on $S$
representing the fixed configuration $\cC$ we construct new
discrete set $V_\cC(S)\subset\R{2}$ and new functional $f\mapsto
\hat f_\cC$. Theorem of Eskin and Masur and theorem of Veech~\cite{zorich:Veech:Siegel}
presented above are
valid  for   $\hat   f_\cC$. The  corresponding
Siegel--Veech  constant  $c(\cC)$  responsible for the  quadratic
growth rate $N_\cC(S,L)\sim c(\cC)\cdot\pi L^2$ of  the number of
collections of homologous  saddle connections of the type $\cC$ on an
individual generic flat  surface  $S$ coincides with the constant
$C(\cC)$        in         the        expression        analogous
to~\eqref{zorich:eq:Veech:counting:theorem}.

Formula~\eqref{zorich:eq:Veech:counting:theorem}  can  be applied
to $\hat{\chi}_L$ for any value of $L$. In particular, instead of
taking large $L$ we can choose a very  small $L=\varepsilon
\ll  1$.  The corresponding  function $\hat{\chi}_\varepsilon(S)$
counts  how  many  collections  of  parallel  $\varepsilon$-short
saddle connections (closed geodesics) of  the  type  $\cC$ we can
find on a flat surface $S\in\cH^{comp}_1(\beta)$.
For   the   flat    surfaces    $S$   outside   of   the   subset
$\cH^{\varepsilon}_1(\cC)\subset\cH^{comp}_1(\beta)$ there are no
such     saddle    connections     (closed     geodesics),     so
$\hat{\chi}_\varepsilon(S)=0$. For surfaces $S$  from  the subset
$\cH^{\varepsilon,thick}_1(\cC)$ there is  exactly one collection
like  this,   $\hat{\chi}_\varepsilon(S)=1$.  Finally,  for   the
surfaces    from    the    remaining    (very    small)    subset
$\cH^{\varepsilon,thin}_1(\cC)=
\cH^{\varepsilon}_1(\cC)\setminus\cH^{\varepsilon,thick}_1(\cC)$
one  has  $\hat{\chi}_\varepsilon(S)\ge 1$. Eskin  and
Masur  have   proved  in~\cite{zorich:Eskin:Masur}  that   though
$\hat{\chi}_\varepsilon(S)$ might be  large on $\cH^{\varepsilon,
thin}_1$ the measure of this subset is so small
(see~\cite{zorich:Masur:Smillie}) that
$$
\int_{\cH^{\varepsilon, thin}_1(\cC)}
\hat{\chi}_\varepsilon(S)\ d\nu_1 = o(\varepsilon^2)
$$
and hence
$$
\int_{\cH^{comp}_1(\beta)} \hat{\chi}_\varepsilon(S)\ d\nu_1 =
\Vol\cH^{\varepsilon, thick}_1(\cC)\  +\ o(\varepsilon^2).
$$
This latter volume is almost the same as the volume
$\Vol\cH^\varepsilon_1(\cC)$ of the neighborhood of the cusp $\cC$
evaluated in equation~\eqref{zorich:eq:volume:of:a:cusp} above,
namely,
$\Vol\cH^{\varepsilon, thick}_1(\cC)=
\Vol\cH^\varepsilon_1(\cC)+ o(\varepsilon^2)$ (see~\cite{zorich:Masur:Smillie}).
Taking into consideration that
$$
\int_{\R{2}} \chi_\varepsilon(x,y)\, dx\, dy=\pi\varepsilon^2
$$
and applying Siegel--Veech
formula~\eqref{zorich:eq:Veech:counting:theorem} to
$\chi_\varepsilon$ we finally get
$$
\cfrac{\Vol\cH^{\varepsilon}_1(\cC)}
{\Vol\cH^{comp}_1(d_1,\dots,d_\noz)}
+o(\varepsilon^2)=
\SVc\cdot\pi\varepsilon^2
$$
which  implies   the  following  formula  for  the  Siegel--Veech
constant $\SVc$:

\begin{multline*}
\SVc=\lim_{\varepsilon\to 0} \cfrac{1}{\pi\varepsilon^2}\,
\cfrac{\Vol(\text{``$\varepsilon$-neighborhood of the cusp $\cC$
''})}{\Vol\cH^{comp}_1(\beta)}
=\\=
(\text{explicit combinatorial factor})\cdot
\cfrac{\prod_{j=1}^k\Vol\cH_1(\beta'_k)}
{\Vol\cH^{comp}_1(\beta)}
\end{multline*}

Sums of the Lyapunov exponents $\nu_1+\dots+\nu_g$
discussed in section~\ref{zorich:s:Generic:Geodesics} are closely
related to the Siegel--Veech constants.

\section{Ergodic Components of the Teichm\"uller Flow}

According          to           the          theorems          of
H.~Masur~\cite{zorich:Masur:Annals:82}           and           of
W.~Veech~\cite{zorich:Veech:Annals:82}   Teichm\"uller   geodesic
flow is ergodic on every connected component of  every stratum of
flat surfaces.  Thus,  the  Lyapunov  exponents  $1+\nu_j$ of the
Teichm\"uller  geodesic  flow  responsible   for   the  deviation
spectrum  of   generic   geodesics   on   a   flat  surface  (see
section~\ref{zorich:s:Generic:Geodesics}),    or    Siegel--Veech
constants responsible for  counting of closed geodesics on a flat
surface                                                      (see
section~\ref{zorich:s:Closed:Geodesics:on:Flat:Surfaces})     are
specific for each connected  component  of each stratum. The fact
that  the  strata  $\cH_1(d_1,\dots,d_\noz)$ are not  necessarily
connected was observed by W.~Veech.

In  order  to formulate the classification theorem for  connected
components of  the  strata  $\cH(d_1,\dots,d_\noz)$  we  need  to
describe the classifying  invariants. There are two of them: spin
structure and  hyperellipticity. Both notions are applicable only
to  part  of  the  strata:  flat   surfaces   from   the   strata
$\cH(2d_1,\dots,2d_\noz)$ have  even  or  odd spin structure. The
strata   $\cH(2g-2)$   and   $\cH(g-1,g-1)$   have   a    special
hyperelliptic connected component.

The results of this section are based on the joint work with
M.~Kontsevich~\cite{zorich:Kontsevich:Zorich}.

\paragraph{Spin structure.}

Consider    a     flat     surface    $S$    from    a    stratum
$\cH(2d_1,\dots,2d_\noz)$. Let  $\rho:  S^1  \to  S$  be a smooth
closed path on $S$; here $S^1$ is a standard circle. Note that at
any point of the surfaces $S$ we know where is the ``direction to
the North''. Hence, at any point $\rho(t)=x\in S$ we can  apply a
compass and measure the direction of the tangent vector $\dot x$.
Moving along  our path $\rho(t)$  we make the tangent vector turn
in the  compass. Thus we get a map  $G(\rho):S^1\to S^1$ from the
parameter circle to the circumference of the compass. This map is
called  the   \emph{Gauss   map}.   We  define  the  \emph{index}
$\ind(\rho)$ of the path $\rho$ as a degree  of the corresponding
Gauss map (or, in other words as the algebraic number of turns of
the tangent vector around the compass) taken modulo $2$.
$$
\ind(\rho)= \deg G(\rho) \mod 2
$$

It  is  easy  to  see  that  $\ind(\rho)$  does   not  depend  on
parameterization.  Moreover,  it  does  not  change  under  small
deformations of the  path. Deforming the path more drastically we
may change its  position with respect to conical singularities of
the flat  metric. Say,  the initial path might go  on the left of
$P_k$ and its deformation might pass on the right of  $P_k$. This
deformation changes  the  $\deg  G(\rho)$.  However,  if the cone
angle  at  $P_k$  is  of  the  type  $2\pi(2d_k+1)$,  then  $\deg
G(\rho)\mod 2$ does not change!  This  observation  explains  why
$\ind(\rho)$ is well-defined  for  a free homotopy class $[\rho]$
when  $S\in\cH(2d_1,\dots,2d_\noz)$  (and hence,  when  all  cone
angles are odd multiples of $2\pi$).

Consider a collection of closed  smooth  paths  $a_1, b_1, \dots,
a_g,  b_g$   representing   a   symplectic   basis   of  homology
$H_1(S,\Z{}/2\Z{})$.   We   define  the   \emph{parity   of   the
spin-structure} of  a flat surface  $S\in\cH(2d_1,\dots,2d_\noz)$
as
$$
\phi(S)=\sum_{i=1}^g
\left(\ind(a_i)+1\right)\left(\ind(b_i)+1\right) \mod 2
$$

\begin{NNLemma}
The value  $\phi(S)$ does not  depend on symplectic
basis  of  cycles   $\{a_i,b_i\}$.   It  does  not  change  under
continuous deformations of $S$ in $\cH(2d_1,\dots,2d_\noz)$.
\end{NNLemma}

The lemma above shows  that  the parity of the  spin  structure is an
invariant of  connected  components  of the strata  of those Abelian
differentials (equivalently, flat surfaces), which have zeroes of even degrees
(equivalently, conical points with cone angles which are odd multiples
of $2\pi$).

\paragraph{Hyperellipticity.}

A flat surface $S$ may  have  a symmetry; one specific family  of
such flat surfaces, which  are  ``more symmetric than others'' is
of a  special interest for us. Recall that  there is a one-to-one
correspondence between  flat  surfaces and pairs (Riemann surface
$M$  ,  holomorphic  1-form  $\omega$).  When  the  corresponding
Riemann  surface  is  hyperelliptic the hyperelliptic  involution
$\tau:M\to  M$  acts  on  any  holomorphic   1-form  $\omega$  as
$\tau^\ast\omega=-\omega$.

We say that a flat surface $S$ is a hyperelliptic flat surface if
there  is  an isometry  $\tau:S\to  S$  such  that  $\tau$  is an
involution,   $\tau\circ\tau=\id$,   and  the   quotient  surface
$S/\tau$  is   a   topological   sphere.   In   flat  coordinates
differential of such involution obviously satisfies $D\tau=-\Id$.

In  a   general  stratum  $\cH(d_1,\dots,d_\noz)$
hyperelliptic flat surfaces form  a small subspace  of  nontrivial  codimension.
However, there are two special  strata,  namely,  $\cH(2g-2)$  and
$\cH(g-1,g-1)$,  for  which hyperelliptic flat surfaces  form  entire
hyperelliptic   connected    components   $\cH^{hyp}(2g-2)$   and
$\cH^{hyp}(g-1,g-1)$ correspondingly.

\begin{NNRemark}
Note that in  the  stratum $\cH(g-1,g-1)$ there are hyperelliptic
flat surfaces of two different types.  A hyperelliptic involution
$\tau S\to S$  may fix the  conical points or  might  interchange
them.  It  is not difficult to  show  that for flat surfaces from  the
connected component $\cH^{hyp}(g-1,g-1)$ the hyperelliptic
involution interchanges the conical singularities.

The remaining  family  of  those  hyperelliptic  flat surfaces in
$\cH(g-1,g-1)$, for  which the hyperelliptic involution keeps the
saddle points fixed,  forms  a subspace of nontrivial codimension
in  the  complement $\cH(g-1,g-1)\setminus\cH^{hyp}(g-1,g-1)$.  Thus, the
hyperelliptic connected  component $\cH^{hyp}(g-1,g-1)$ does  not
coincide with the space of all hyperelliptic flat surfaces.
\end{NNRemark}

\paragraph{Classification theorem for Abelian differentials.}

Now, having introduced the classifying invariants  we can present
the classification of  connected  components of strata of
flat surfaces (equivalently, of strata of Abelian
differentials).

\begin{Theorem}[M.~Kontsevich and A.~Zorich]
All connected components of any stratum  of flat surfaces
of genus $g\ge 4$ are described by the following list:

The stratum $\mathcal{H}(2g-2)$ has  three  connected components:
the   hyperelliptic   one,  $\mathcal{H}^{hyp}(2g-2)$,   and  two
nonhyperelliptic   components:   $\mathcal{H}^{even}(2g-2)$   and
$\mathcal{H}^{odd}(2g-2)$  corresponding  to even  and  odd  spin
structures.

The stratum  $\mathcal{H}(2d,2d)$,  $d\ge  2$ has three connected
components:  the  hyperelliptic one,  $\mathcal{H}^{hyp}(2d,2d)$,
and two  nonhyperelliptic components: $\mathcal{H}^{even}(2d,2d)$
and $\mathcal{H}^{odd}(2d,2d)$.

All      the       other       strata       of      the      form
$\mathcal{H}(2d_1,\dots,2d_\noz)$ have  two connected components:
$\mathcal{H}^{even}(2d_1,\dots,2d_\noz)$                      and
$\mathcal{H}^{odd}(2d_1,\dots,2d_n)$, corresponding to  even  and
odd spin structures.

The stratum $\mathcal{H}(2d-1,2d-1)$, $d\ge 2$, has two connected
components;  one   of  them:  $\mathcal{H}^{hyp}(2d-1,2d-1)$   is
hyperelliptic;  the  other  $\mathcal{H}^{nonhyp}(2d-1,2d-1)$  is
not.

All the other strata of flat surfaces
of genera $g\ge 4$ are nonempty and connected.
\end{Theorem}

In the case of small genera $1\le g\le 3$ some components are
missing in comparison with the general case.

\addtocounter{Theorem}{-1}
\setcounter{TheoremPrime}{\theTheorem}
\begin{TheoremPrime}
The moduli space  of flat surfaces  of  genus
$g=2$    contains    two     strata:    $\mathcal{H}(1,1)$    and
$\mathcal{H}(2)$.  Each  of them is connected and coincides  with
its hyperelliptic component.

Each of  the  strata  $\mathcal{H}(2,2)$, $\mathcal{H}(4)$ of the
moduli space  of flat surfaces of genus $g=3$
has two  connected  components:  the  hyperelliptic  one, and one
having odd spin structure. The other  strata  are  connected  for
genus $g=3$.
\end{TheoremPrime}

Since there  is  a  one-to-one  correspondence  between connected
components  of  the  strata  and  extended   Rauzy  classes,  the
classification theorem  above  classifies also the extended Rauzy
classes.

Connected  components  of  the strata $\cQ(d_1,\dots,d_\noz)$  of
meromorphic quadratic differentials with at most simple poles are
classified in the paper of E.~Lanneau~\cite{zorich:Lanneau}.

\medskip
\paragraph{Bibliographical notes.}

As a much  more  serious accessible introduction to Teichm\"uller
dynamics   I   can   recommend   a  collection  of   surveys   of
A.~Eskin~\cite{zorich:Eskin:Handbook},
G.~Forni~\cite{zorich:Forni:Handbook},       P.~Hubert        and
T.~Schmidt~\cite{zorich:Hubert:Schmidt:Handbook}              and
H.~Masur~\cite{zorich:Masur:Handbook:1B}, organized  as a chapter
of the Handbook of Dynamical Systems.  I  also  recommend  recent
surveys               of               H.~Masur               and
S.~Tabachnikov~\cite{zorich:Masur:Tabachnikov}       and       of
J.~Smillie~\cite{zorich:Smillie:billiards}   especially   in  the
aspects related to billiards in  polygons.  The  part  concerning
renormalization   and  interval   exchange   transformations   is
presented           in           the          survey           of
J.-C.~Yoccoz~\cite{zorich:Yoccoz:Les:Houches}.     The      ideas
presented in the current  paper  are illustrated in more detailed
way in the survey~\cite{zorich:Zorich:Les:Houches}.

\medskip
\paragraph{Acknowledgements.}

Considerable part of results presented in this survey is obtained in
collaboration. I use this opportunity to thank A.~Eskin, M.~Kontsevich
and H.~Masur for the pleasure to work with them.
I am grateful to M.-C.~Vergne for her help in preparation of
pictures. I highly  appreciate careful and
responsible editing by I.~Zimmermann.


\end{document}